 \documentclass[12pt]{article}
 \usepackage{amsmath,amsthm,amsfonts,amssymb}
\usepackage{verbatim}
\usepackage{latexsym}
\usepackage[mathscr]{euscript}
\newtheorem{theorem}{Theorem}[section]
\newtheorem{lemma}[theorem]{Lemma}
\newtheorem{proposition}[theorem]{Proposition}
\newtheorem{definition}[theorem]{Definition}
\newtheorem{corollary}[theorem]{Corollary}

\newtheorem{remark}[theorem]{Remark}
\newenvironment{prf} {{\bf Proof.}}{\hfill $\Box$}

\newcommand\RR{{\mathbb{R}}}

\newcommand\NN{\mathbb{ N}}
\newcommand\ZZ{\mathbb{ Z}}

\newcommand\HH{\mathbb{ H}}

\newcommand{\note}[2][\null]{%
  \marginpar{\renewcommand{\baselinestretch}{1}\vspace{-1em}\hrule\vspace{3pt}%
  \tiny\raggedright{#2\ifx#1\null\else\\\hfill---
  {\em #1}\fi}\vspace{1.5em}}%
}

\def\sideremark#1{\ifvmode\leavevmode\fi\vadjust{\vbox to0pt{\vss
 \hbox to 0pt{\hskip\hsize\hskip1em
\vbox{\hsize2cm\tiny\raggedright\pretolerance10000 
 \noindent #1\hfill}\hss}\vbox to8pt{\vfil}\vss}}} 

\begin{document}
\title{ Homogeneous Besov spaces on stratified Lie groups and their wavelet characterization}
\author{
Hartmut F\"uhr\\
\footnotesize\texttt{{fuehr@matha.rwth-aachen.de}}\\
  Azita Mayeli \\
\footnotesize\texttt{{amayeli@qcc.cuny.edu}} }
 \maketitle

\begin{abstract}
We establish wavelet characterizations of homogeneous Besov spaces on stratified Lie groups, both in terms of continuous and discrete wavelet systems.

We first introduce a notion of homogeneous Besov space $\dot{B}_{p,q}^s$ in terms of a Littlewood-Paley-type decomposition, in analogy to the well-known characterization of the Euclidean case. Such decompositions can be defined via the spectral measure of a suitably chosen sub-Laplacian. We prove that the scale of Besov spaces  is independent of the precise choice of Littlewood-Paley decomposition. In particular, different sub-Laplacians yield the same Besov spaces.

We then turn to wavelet characterizations, first via continuous wavelet transforms (which can be viewed as continuous-scale Littlewood-Paley decompositions), then via discretely indexed systems. We prove the existence of wavelet frames and associated atomic decomposition formulas for all homogeneous Besov spaces ${\dot B}_{p,q}^{s}$, with $1 \le p,q < \infty$ and $s \in \mathbb{R}$.

 \end{abstract}
  {\footnotesize



 }

\section{Introduction}\label{introduction}

To a large extent, the success of wavelets in applications can be attributed to the realization that wavelet bases are universal unconditional bases for a large class of smoothness spaces, including all homogeneous Besov spaces. Given a wavelet orthonormal basis $\{ \psi_{j,k} \}_{j,k} \subset {\rm L}^2(\mathbb{R}^n)$ (consisting of sufficiently regular wavelets with vanishing moments) and $f \in {\rm L}^2(\mathbb{R}^n)$, the expansion
\[
  f = \sum_{j,k} \langle f, \psi_{j,k} \rangle \psi_{j,ĸk}
\] converges not only in $\| \cdot \|_{L^2}$, but also in any other Besov space norm $\| \cdot \|_{\dot{B}_{p,q}^s}$, as soon as $f$ is contained in that space. Furthermore, the latter condition can be read off the decay behaviour of the wavelet coefficients
$\{\langle f, \psi_{j,k} \rangle \}_{j,k}$ associated to $f$ in a straightforward manner.

This observation provided important background and heuristics for
many wavelet-based methods in applications such as denoising and
data compression, but it was also of considerable theoretical
interest, e.g. for the study of operators. In this paper we provide
similar results for simply connected stratified Lie groups. To our
knowledge, studies of Besov spaces in this context have been
restricted to the inhomogeneous cases: The definition of
inhomogeneous Besov spaces on stratified Lie groups was introduced
independently by Saka \cite{Saka79}, and in a somewhat more general
setting by Pesenson \cite{Pe79,Pe83}. Since then, the study of Besov
spaces on Lie groups remained restricted to the inhomogeneous cases
\cite{Tr1,Tr2,giulini,Skrzypczak02,Furio-Melzi-Veneruso06}, with the
notable exception of \cite{Ba} which studied homogeneous Besov
spaces on the Heisenberg group. A further highly influential source
for the study of function spaces associated to the sub-Laplacian is
Folland's paper \cite{Folland75}.

The first wavelet systems on stratified Lie groups (fulfilling certain technical assumptions) were constructed by Lemari\'e \cite{lemarie}, by suitably adapting concepts from spline theory. Lemari\'e also indicated that the wavelet systems constructed by his approach were indeed unconditional bases of Saka's inhomogeneous Besov spaces. Note in particular that an adaptation, say, of the arguments in \cite{FrazierJawerth85} for a proof of such a characterization requires a sampling theory for bandlimited functions on stratified groups, which was established only a few years ago by Pesenson \cite{Pesenson}; see also \cite{FuGr}.

More recent constructions of both continuous and discrete wavelet systems were based on the spectral theory of the sub-Laplacian \cite{gm1}. Given the central role of the sub-Laplacian both in \cite{Furio-Melzi-Veneruso06} and \cite{gm1}, and in view of Lemari\'e's remarks, it seemed quite natural to expect a wavelet characterization of homogeneous Besov spaces, and it is the aim of this paper to work out the necessary details.

The paper is structured as follows: After  reviewing the basic
notions concerning stratified Lie groups and their associated
sub-Laplacians in Section \ref{preliminary-and-notations}, in
Section \ref{besov-spaces} we introduce a Littlewood-Paley-type
decomposition of functions and tempered discributions on $G$. It is
customary to employ the spectral calculus of a suitable
sub-Laplacian for the definition of such decompositions, see e.g.
\cite{Furio-Melzi-Veneruso06}, and this approach is also used here
(Lemma \ref{lem:construct_LP_adm}). However, this raises the issue
of consistency: The spaces should reflect properties of the group,
not of the sub-Laplacian used for the construction of the
decomposition. Using a somewhat more general notion than the
$\phi$-functions in \cite{FrazierJawerth85} allows to establish that
different choices of sub-Laplacian result in the same scale of Besov
spaces (Theorem \ref{thm:Besov_norm_equiv}). In Section 4, we derive
a characterization of Besov spaces in terms of continuous wavelet
transform, with a wide variety of wavelets to choose from (Theorem
\ref{thm:Char_CWT}). As a special case one obtains a
characterization of homogeneous Besov spaces in terms of the heat
semigroup. (See the remarks before Theorem \ref{thm:Char_CWT}.)

In Section  \ref{Abtast-Teil}, we study discrete characterizations of Besov spaces obtained by sampling the Calder\'{o}n decomposition.  For this purpose, we introduce the coefficient space $\dot{b}_{p,q}^s$. The chief result is Theorem \ref{thm:besov_discrete}, establishing that the wavelet coefficient sequence of $f \in \dot{B}_{p,q}^s$ lies in
$\dot{b}_{p,q}^s$. Section \ref{Abtast-Teil} introduces our most important tool to bridge the gap between continuous and discrete decompositions, namely oscillation estimates.

We then proceed to study wavelet synthesis and frame properties of the wavelet system. Our main result in this respect is that for all sufficiently dense regular sampling sets $\Gamma$, the discrete wavelet system $(\psi_{j,\gamma})_{j \in \mathbb{Z},\gamma \in \Gamma}$ obtained by shifts from $\gamma$ and dilations by powers of $2$ is a {\em universal Banach frame} for all Besov spaces. In other words, the wavelet system allows the decomposition
\[
  f = \sum_{j,\gamma} r_{j,\gamma} \psi_{j,\gamma}
\] converging unconditionally in ${\dot B}_{p,q}^s$ whenever $f \in {\dot B}_{p,q}^s$, with coefficients $\{ r_{j,\gamma} \}_{j,\gamma} \in {\dot b}_{p,q}^s$ depending linearly and boundedly on $f$, and satisfying the norm equivalence
\[
 \| \{ r_{j,\gamma} \}_{j,\gamma} \|_{{\dot b}_{p,q}^s} \asymp \| f \|_{{\dot B}_{p,q}^s}~.
\]

\section{Preliminaries and Notation}\label{preliminary-and-notations}

Following  the terminology in \cite{FollandStein82}, we call a Lie
group ${G}$ stratified if it is  connected and simply
connected, and its Lie algebra $\mathfrak{g}$ decomposes as a direct
sum $\mathfrak{g}= V_1\oplus \cdots \oplus V_m $
$\left[V_1,V_k\right]  = V_{k+1}$ for $1\leq k<m$ and $\left[
V_1,V_m\right]=\{0\}$. Then $\mathfrak{g}$ is nilpotent of
step $m$, and generated as a Lie algebra by $V_1$.  
Euclidean spaces $\RR^n$ and the Heisenberg group $\HH^n$ are
 examples of stratified Lie groups.

If ${G}$ is stratified, its Lie algebra admits a
canonical (natural)  family of dilations, namely
\begin{align}\notag
\delta_r(X_1+X_2+\cdots+ X_m)= rX_1+r^2X_2+\cdots +r^mX_m \quad (X_j\in V_j)~~(r>0)~,\notag
\end{align}
which are Lie algebra automorphisms.
We identify $G$ with $\mathfrak{g}$ through the exponential map.
 Hence $G$ is a Lie group with underlying
manifold $\RR^n$, for some $n$, and the group product provided by the Campbell-Baker-Hausdorff formula.
The dilations are then also group automorphisms of $G$. Instead of writing $\delta_a (x)$ for $x \in G$ and
$a >0$, we simply use $ax$, whenever a confusion with the Lie group product is excluded. After choosing a basis of $\mathfrak{g}$ obtained as a union of bases of the $V_i$, and a possible change of coordinates, one therefore has for $x\in G$ and $a>0$ that
\begin{equation}\label{eq:dil}
ax=(a^{d_1}x_1,\cdots,a^{d_n}x_n),
\end{equation}
for integers  $d_1\leq \cdots\leq d_n$, according to $x_i \in V_{d_i}$. 

Under our identification of $G$ with $\mathfrak{g}$, polynomials on
$G$ are polynomials on $\mathfrak{g}$ (with respect to any linear
coordinate system on the latter). Polynomials on $G$ are written as
\begin{equation} \label{eqn:poly_gen_form} p\left(\sum_{i=1}^{\dim(G)}  x_i Y_i\right) =  \sum_I c_I x^I
 \end{equation}
where $c_I \in \mathbb{C}$ are the coefficients, and $x^I=
x_1^{I_1}x_2^{I_2}\cdots x_n^{I_n}$ the monomials associated to the
multi-indices $I \in \mathbb{N}^{\{ 1,\ldots,n\}}$. For a
multi-index $I$, define
\[
 d(I) = \sum_{i=1}^n I_i n(i)~,n(i) = j \mbox{ for } Y_i \in V_j~.
\] A polynomial of the type (\ref{eqn:poly_gen_form}) is called {\em of homogeneous degree $k$} if $d(I) \le k$ holds, for all multiindices $I$ with $c_I \not= 0$. We write $\mathcal{P}_k$ for the space of polynomials of homogeneous degree $k$.

We let $\mathcal{S}(G)$
 denote the space of Schwartz functions on
$G$. By definition $\mathcal{S}(G) = \mathcal{S}(\mathfrak{g})$. Let
$\mathcal{S}'(G)$ and  $\mathcal{S}'(G)/\mathcal{P}$ denote the
space of distributions and distributions modulo polynomials on $G$,
respectively. The duality between the spaces is denoted by the map
$\left( \cdot, \cdot \right): \mathcal{S}' (G) \times \mathcal{S}(G)
\to \mathbb{C}$. Most of the time, however, we will work with the
sesquilinear version $\langle f, g \rangle = (f, \overline{g})$, for
$f \in \mathcal{S'}(G)$ and $g \in \mathcal{S}(G)$.

Left Haar measure on $G$ is induced by Lebesgue
measure on its Lie algebra, and it is also right-invariant.
 The number $Q= \sum_1^m j(\dim V_j)$ will be called the \textsl{
homogeneous dimension} of $G$.
 (For  instance, for $G=\RR^n$ and $\HH^n$ we have $Q=n$ and $Q=2n+2$, respectively.)
 For any  function $\phi$  on ${G}$ and $a>0$, 
the ${\rm L}^1$-normalized dilation  of $\phi$ is defined by
\begin{align}\notag
D_a \phi (x)= a^{Q}\phi(ax).
\end{align}
Observe that this action preserves the  $L^1$-norm, i.e., $\parallel
D_a \phi \parallel_1= \parallel \phi\parallel$.  We fix a
homogeneous quasi-norm $|\cdot |$ on $G$ which is smooth away from
$0$,  $|ax| = a|x|$ for all $x \in G$, $a \geq 0$, $|x^{-1}| = |x|$
for all $x \in G$, with $|x|
> 0$ if $x \neq 0$, and fulfilling a triangle inequality $|xy| \le C(|x|+|y|)$, with constant $C>0$.
Confer \cite{FollandStein82} for the construction of homogeneous norms, as well as further properties. 

Moreover, by Proposition 1.15 \cite{FollandStein82}, for any $r> 0$,
there is a finite $C_r > 0$ such that
$\int_{|x| > R} |x|^{-Q-r} dx= C_r R^{-r}$ for all $R > 0$.

 Our conventions for left-invariant operators on $G$ are as follows: We let $Y_1,\ldots,Y_n$ denote a basis of
 $\mathfrak{g}$, obtained as a union of bases of the $V_i$. In particular, $Y_1,\ldots,Y_l$, for $l = \dim (V_1)$, is a basis of $V_1$. Elements of the Lie algebra are identified in the usual manner with left-invariant differential operators on $G$. Given a multi-index $I\in \mathbb{N}_0^n$, we write $Y^I$ for $Y_1^{I_1} \circ \ldots \circ Y_n^{I_n}$. A convenient characterization of Schwartz functions in terms of left-invariant operators states that $f \in \mathcal{S}(G)$ iff, for all $N \in \mathbb{N}$, $|f|_N < \infty$, where
\[
 |f|_N = \sup_{|I| \le N, x \in G} (1+|x|)^{N} |Y^I f(x)|~.
\] In addition, the norms $|\cdot |_N$ induce the topology of $\mathcal{S}(G)$ (see \cite{FollandStein82}).

The sub-Laplacian operator on $G$ can be viewed as
the analog of the Laplacian operator on
  $\RR^n$ defined by $L=-\sum_{i=1}^n
 \frac{\partial^2}{\partial x_k^2}$. Using the above conventions for the choice of basis $Y_1,\ldots,Y_n$ and $l = \dim(V_1)$, the sub-Laplacian is defined as $L=-\sum_{i=1}^l Y_i^2$.
Note that a less restrictive notion of sub-Laplacians can also be found in the literature (e.g., any sum of squares of Lie algebra generators); we stress that the results in this paper crucially rely on the definition presented here.
A linear differential operator $T$ on $G$ is called homogenous of
degree $l$ if $T(f\circ \delta_a)= a^l (Tf)\circ \delta_a$ for any
$f$ on $G$. By choice of the $Y_i$ for $i \le l$, these operators
are homogeneous of degree 1; it follows that $L$ is  homogenous of
degree $2$, and $L^k$ is homogenous of degree $2k$. Furthermore, any
operator of the form $Y^I$ is homogeneous of degree $d(I)$.

 When restricted to $C_c^{\infty}$,
$L$ is formally self-adjoint: for any $f,g\in C_c^\infty(G)$, $\langle
 Lf, g\rangle= \langle f, Lg\rangle$. (For more see \cite{gm1}).  Its closure has domain
$\mathcal{D} = \{u \in L^2(G): Lu \in L^2(G)\}$, where we take
$Lu$ in the sense of distributions.  From this fact it quickly
follows that this closure is self-adjoint and is in fact the unique
self-adjoint extension of $L|_{C_c^{\infty}}$; we denote this extension also by the symbol $L$.

Suppose that $L$ has spectral resolution
\begin{align}\notag
L= \int_0^\infty \lambda dP_\lambda
\end{align}
where $dP_\lambda$ is the  projection measure.
 For  a bounded Borel function  $\hat f$  on $[0,\infty)$,
the operator
\begin{align}\label{integral} \hat f(L)=\int_0^\infty
\hat f(\lambda)dP_\lambda
\end{align}
is a bounded integral operator on $L^2({G})$ with a convolution
distribution kernel in $L^2(G)$ denoted by $f$, and
\begin{align}\notag
\hat f(L)\eta= \eta\ast  f \quad \forall \;\eta\in \mathcal{S}(G)~.
\end{align}
An important fact to be used later on is that for rapidly decaying
smooth functions, $f \in \mathcal{S}(\mathbb{R}^+)$, the kernel
associated to $\hat{f}(L)$ is a Schwartz function. For a function
$f$ on $G$ we define $\tilde{f}(x) = f(x^{-1})$ and $f^* =
\overline{\tilde{f}}$. We will frequently use that for $f \in {\rm
L}^2(G) \cap {\rm L}^1(G)$, the adjoint of the convolution operator
$g \mapsto g \ast f$ is provided by $g \mapsto g \ast f^*$.

%
%
%
%


\section{Homogeneous Besov Spaces on Stratified Lie Groups 
}\label{besov-spaces}
In this section we define homogeneous Besov spaces on
stratified Lie groups via Littlewood-Paley decompositions of distributions $u$ as
\[
 u = \sum_{j \in \mathbb{Z}} u \ast \psi_j^* \ast \psi_j~,
\] where $\psi_j$ is a dilated copy of a suitably chosen Schwartz function
$\psi$. In the Euclidean setting, it is customary to construct
$\psi$ by picking a dyadic partition of unity on the Fourier
transform side and applying Fourier inversion. The standard way of
transferring this construction to stratified Lie groups consists in
replacing the Fourier transform by the spectral decomposition of a
sub-Laplacian $L$, see Lemma \ref{lem:construct_LP_adm} below.
However, this approach raises the question to which extent the
construction depends on the choice of $L$. It turns out that the
precise choice of sub-Laplacian obtained from a basis of $V_1$ is irrelevant. 
In order to prove this, we study Littlewood-Paley
decompositions in somewhat different terms. The right setting for
the study of such decompositions is the space of tempered
distributions modulo polynomials, and the easiest approach to this
convergence is via duality to a suitable space of Schwartz
functions.

\begin{definition}
 Let $N \in \mathbb{N}$. A function $f: G \to \mathbb{C}$ has polynomial
  decay order $N$ if  there exists a constant $C>0$ such
 that, for all $x \in G$
\[
 |f(x)| \le C (1+|x|)^{-N}
\]
$f$ has vanishing moments of order $N$, if  one has
\begin{equation}
 \label{eqn:defn_van_mom}\forall p \in \mathcal{P}_{N-1}~:~\int_G f(x) p(x) dx = 0~,
\end{equation} with absolute convergence of the integral.
\end{definition}
Under our identification of $G$ with $\mathfrak{g}$, the inversion map $x \mapsto x^{-1}$ is identical to the additive inversion map. I.e., $x^{-1} = -x$, and it follows that $\tilde{p} \in \mathcal{P}_{N}$ for all $p \in \mathcal{P}_N$. Thus, if $f$ has vanishing moments of order $N$, then for all $p \in \mathcal{P}_{N-1}$
\[
  \int_G \tilde{f}(x) p(x) dx = \int_G f(x) \tilde{p}(x) dx  = 0,
\] i.e., $\tilde{f}$ has vanishing moments of order $N$ as well.

Vanishing moments are central to most estimates in wavelet analysis, by the following principle: In a convolution product of the type $g \ast D_t f$, vanishing moments of one factor together with smoothness of the other result in decay. Later on, we will apply the lemma to Schwartz functions $f,g$, where only the vanishing moment assumptions are nontrivial. 
The more general version given here is included for reference.

\begin{lemma} \label{lem:decay_conv} Let $N,k \in \mathbb{N}$ be arbitrary.
 \begin{enumerate}
\item[(a)]  Let $f \in C^{k}$, such that $Y^I(f)$ is of decay order $N$, for all $I$ with $d(I) \le k$. Let $g$ have vanishing moments of order $k$ and decay order $N+k +Q+1$. Then there exists a constant, depending only on the decay of $Y^I(f)$ and $g$, such that
\begin{equation} \label{eqn:decay_conv1}
 \forall x \in G~ \forall~ 0 < t < 1 ~:~|g \ast (D_t f)(x)| \le C t^{k+Q} (1+|tx|)^{-N}~.
\end{equation}
 In particular, if $p>Q/N$,
 \begin{equation} \label{eqn:decay_conv1_norm}
 \forall x \in G~ \forall~ 0 < t < 1 ~:~\| g \ast (D_t f) \|_p \le C' t^{k+Q(1-1/p)} ~.
\end{equation}
 \item[(b)] Now suppose that $g \in C^k$, with $Y^I(\tilde{g})$ of  decay order $N$ for all $I$ with $d(I) \le k$. Let $f$ have vanishing moments of order $k$ and decay order $N+k+Q+1$. Then there exists a constant, depending only on the decay of $f$ and $Y^I(\tilde{g})$, such that
\begin{equation} \label{eqn:decay_conv2}
  \forall x \in G~ \forall~ 1 < t < \infty ~:~|g \ast (D_t f)(x)| \le C t^{-k}(1+|x|)^{-N} ~.
\end{equation}
 In particular, if $p>Q/N$,
 \begin{equation} \label{eqn:decay_conv2_norm}
 \forall x \in G~ \forall~ 1 < t < \infty ~:~\|g \ast (D_t f)\|_p \le C' t^{-k} ~.
\end{equation}
\end{enumerate}
\end{lemma}
\begin{proof}
First, let us prove $(a)$: Let $0<t<1$.
For $x \in G$, let $P_{x,D_t \tilde{f}}^k$ denote the left Taylor polynomial of $D_t\tilde{f}$ with homogeneous degree $k-1$, see \cite[1.44]{FollandStein82}. By that result,
\[
 \left|f(y^{-1}x) - P^k_{x,D_t  \tilde{f}}(y)\right| \le C_k |y|^{k} \sup_{|z| \le b^{k} |y|,d(I) = k} |Y^I (D_t \tilde{f})(xz)|~,
\] with suitable positive constants $C_k$ and $b$. We next use the homogoneity properties of the partial derivatives\cite[p.21]{FollandStein82}, together with the decay condition on $Y^I f$ to estimate for $I$ with $d(I)=k$
\begin{eqnarray*}
\sup_{|z| \le b^{k} |y|}  |Y^I (D_t \tilde{f})(xz)| & = & t^{k} \sup_{|z| \le b^k |y|} |D_t (Y^I \tilde{f})(xz)| \\
 & = & t^{k+Q} \sup_{|z| \le b^k y} |(Y^I \tilde{f})(t(x\cdot z))| \\
 & \le &  t^{k+Q} \sup_{z \le b^{k}|y|} C_f ( 1+|t(x \cdot z)|)^{-N} \\
 & \le &  t^{k+Q} \sup_{z \le b^{k}|y|} C_f (1+|tx|)^{-N} (1+|tz|)^N \\
 & \le &   t^{k+Q} (1+b)^{kN} C_f (1+|tx|)^{-N} (1+ |y|)^N~,
\end{eqnarray*} where the penultimate inequality used \cite[1.10]{FollandStein82}, and the final estimate used $|ty| = t |y| \le |y|$.  Thus
\[
  \left|f(y^{-1}x) - P^k_{x,D_t\tilde{f}}(y)\right| \le \tilde{C}_k  t^{k+Q} (1+|y|)^{N+k} (1+|tx|)^{-N}~.
\] Next, using vanishing moments of $g$,
\begin{eqnarray*}
\lefteqn{\left| (g \ast D_t f) (x) \right|} \\
 & \le & \int_G |g(y)| \left| D_t f(y^{-1}x) -P_{x,D_t \tilde{f}}^k(y)\right| dy \\
 & \le & \tilde{C}_k (1+|tx|)^{-N}  t^{k+Q} \int_G |g(y)| ~ (1+|y|)^{N+k} dy \\
& \le & \tilde{C}_k (1+|tx|)^{-N}   t^{k+Q}  \int_G C_g (1+|y|)^{-Q-1} dy~,
\end{eqnarray*} and the integral is finite by \cite[1.15]{FollandStein82}. This proves (\ref{eqn:decay_conv1}), and
(\ref{eqn:decay_conv1_norm}) follows by
\[
 \| g \ast D_t f\|_p \le C' t^{k+Q} \left( \int_G (1+|tx|)^{-Np} dx \right)^{1/p} \le C'' t^{k+Q-Q/p}  ~,
\] using $Np > Q$.

For part (b), we first observe that
\[
  (g \ast D_t f)(x) = t^Q \left(\tilde{f} \ast D_{t^{-1}} \tilde{g} \right) (t.x)~.
\] Our assumptions on $f,g$ allow to invoke part (a) with $\tilde{g},\tilde{f}$ replacing $f,g$, and (\ref{eqn:decay_conv2}) follows immediately. (\ref{eqn:decay_conv2_norm}) is obtained from this by straightforward integration.
\end{proof}

We let $\mathcal{Z}(G)$ denote the space of Schwartz functions with all moments vanishing.
  We next consider properties of $\mathcal{Z}(G)$ as a subspace of
$\mathcal{S}(G)$ with the relative topology.

\begin{lemma}\label{identification-dual-spaces} $\mathcal{Z}(G)$ is a closed subspace (in particular complete) of $\mathcal{S}(G)$, with $\mathcal{S} (G) \ast \mathcal{Z}(G) \subset \mathcal{Z}(G)$, as well as $\tilde{f} \in \mathcal{Z}(G)$ for all $f \in \mathcal{Z}(G)$. The topological dual of  $\mathcal{Z}(G)$,  $\mathcal{Z}'(G)$,
can be canonically identified with the factor space $\mathcal{S}'(G)/\mathcal{P}$.
\end{lemma}
\begin{proof}
By definition, $\mathcal{Z}(G)$ is the intersection of kernels of a
family of tempered distributions, hence a closed subspace. For $p
\in \mathcal{P}$ and $f \in \mathcal{Z}(G)$, one has by
unimodularity of $G$ that $\langle p,\tilde{f}\rangle = \langle
\tilde{p},f \rangle = 0$, since $\tilde{p}$ is a polynomial. But
then, for any $g \in \mathcal{S}(G)$ and $f \in \mathcal{Z}(G)$, one
has for all polynomials $p$ on $G$ that
\[
 \langle g \ast f, p \rangle = \langle g , p \ast \tilde{f} \rangle = \langle g, 0 \rangle  = 0~,
\] since $f \in \mathcal{Z}(G)$ implies $p \ast \tilde{f} = 0$ (translation on $G$ is polynomial). Thus $g \ast f \in \mathcal{Z}(G)$.

All further properties of $\mathcal{Z}(G)$ follow from the
corresponding statements concerning $\mathcal{Z}(\mathbb{R}^n)$. For
identification of $\mathcal{Z}'(\mathbb{R}^n)$ with the quotient
space $\mathcal{S}(\mathbb{R}^n)'/\mathcal{P}$,  we first observe
that a tempered distribution $\varphi$ vanishes on
$\mathcal{Z}(\mathbb{R}^n)$ iff its (Euclidean) Fourier transform is
supported in $\{ 0 \}$, which is well-known to be the case iff
$\varphi$ is a polynomial. Using this observation, we map $u \in
\mathcal{Z}'(\mathbb{R}^n)$ to $\tilde u+ \mathcal{P}$, where
$\tilde u $ is a continuous extension of $u$ to all of
$\mathcal{S}(\mathbb{R}^n)$; such an extension exists by the
Hahn-Banach theorem. The map is well-defined because the difference
of two extensions of $u$ annihilate $\mathcal{Z}(\mathbb{R}^n)$,
hence differ by a polynomial. Linearity follows from
well-definedness. Furthermore, the inverse of the mapping is clearly
obtained by assigning $w + \mathcal{P}$ to the restriction
$w|_{\mathcal{Z}(G)}$.
%
\end{proof}
In the following, we will usually not explicitly distinguish between
$u \in \mathcal{S}'(G)$
 and its equivalence class modulo polynomials,
and we will occasionally write $u \in \mathcal{S}'(G)/\mathcal{P}$.
The topology of $\mathcal{S}'(G)/\mathcal{P}$ is just the topology
of pointwise convergence on the elements of $\mathcal{Z}(G)$: For
any net $(u_j + \mathcal{P} )_{j \in I}$, $u_j + \mathcal{P} \to u +
\mathcal{P}$ holds if and only  if $\langle u_j, \varphi \rangle \to
\langle u, \varphi \rangle$, for all $\varphi \in\mathcal{Z}(G)$. We next study
convolution on $\mathcal{S}'(G)/\mathcal{P}$.

\begin{lemma} \label{lem:conv_SmP}
  For every $\psi \in \mathcal{S}(G)$, the map $u \mapsto u \ast \psi$ is a well-defined
  and continuous operator $\mathcal{S}'(G)/\mathcal{P} \to \mathcal{S}'(G)/\mathcal{P}$.
   If $\psi \in \mathcal{Z}(G)$, the associated convolution operator is a well-defined and continuous operator $\mathcal{S}'(G)/\mathcal{P} \to \mathcal{S}'(G)$.
\end{lemma}
\begin{proof}
Note that $\mathcal{P} \ast \mathcal{S}(G) \subset \mathcal{P}$.
Hence $u \mapsto u \ast \psi$ induces a well-defined canonical map
$\mathcal{S}'(G)/\mathcal{P} \to \mathcal{S}'(G)/\mathcal{P}$.
Furthermore, $u \mapsto u \ast \psi$ is continuous on
$\mathcal{S}'(G)$, as a consequence of \cite[Proposition
1.47]{FollandStein82}. Therefore, for any net $u_j \to u$ and any
$\varphi \in \mathcal{Z}(G)$, the fact that $\varphi \ast \psi^* \in
\mathcal{Z}(G)$ allows to write
\begin{equation} \label{eqn:conv_cont} \langle u_j \ast \psi,\varphi \rangle = \langle u_j ,\varphi \ast \psi^* \rangle \to \langle u, \varphi \ast \psi^* \rangle =
 \langle u \ast \psi, \varphi \rangle~,
\end{equation} showing $u_j \ast \psi \to u \ast \psi$ in $\mathcal{S}'(G)/\mathcal{P}$.

For $\psi \in \mathcal{Z}(G)$, the fact that $\mathcal{P} \ast \psi
= \{ 0 \}$ makes the mapping $u \mapsto u \ast \psi \in
\mathcal{S}'(G)$ well-defined modulo polynomials. The continuity
statement is proved by (\ref{eqn:conv_cont}), with assumptions on
$\psi$ and $\varphi$ switched.
\end{proof}

The definition of homogeneous Besov spaces requires taking ${\rm L}^p$-norms of elements of $\mathcal{S}'(G)/\mathcal{P}$. The following remark clarifies this.

\begin{remark} \label{rem:LP_norm_dist}
Throughout this paper, we use the canonical embedding  $L^p(G) \subset \mathcal{S}'(G)$. For $p< \infty$ this gives rise to an embedding $L^p(G) \subset \mathcal{S}'(G)/\mathcal{P}$, using that $\mathcal{P} \cap {\rm L}^p(G) = \{ 0 \}$. Consequently, given $u \in \mathcal{S}'(G)/\mathcal{P}$, we let
\begin{equation} \label{eqn:defn_lp_P} \| u \|_p = \| u + q \|_p \mbox{ whenever }u+q \in {\rm L}^p(G)~, \mbox{ for suitable } q \in \mathcal{P}
\end{equation} assigning the value $\infty$ otherwise.
Here the fact that $\mathcal{P} \cap {\rm L}^p(G) = \{ 0 \}$ guarantees that the decomposition is unique, and thus (\ref{eqn:defn_lp_P}) well-defined.

By contrast, $\| \cdot \|_\infty$ can only be defined on $\mathcal{S}'(G)$, if we assign the value $\infty$ to $u \in \mathcal{S}'(G) \setminus {\rm L}^\infty(G)$.

Note that with these definitions,  the Hausdorff-Young inequality $\| u \ast f \|_p \le \|u \|_p \| f \|_1$ remains valid for all $f \in \mathcal{S}(G)$, and all  $u \in \mathcal{S}'(G)/\mathcal{P}$ (for $p<\infty$), resp. $u \in \mathcal{S}'(G)$ (for $p=\infty$).
For $p=\infty$, this is clear. For $p< \infty$, note that if $u + q \in {\rm L}^p(G)$, then $(u+q) \ast \psi = u \ast \psi + q \ast \psi \in {\rm L}^p(G)$ with $q \ast \psi \in \mathcal{P}$.
\end{remark}

We now introduce a general Littlewood-Paley-type decomposition.
For this purpose we define for $\psi \in \mathcal{S}(G)$:
\[
 \psi_j = D_{2^j} \psi ~.
\]
\begin{definition} \label{defn:LP_S(G)}
 A  function $\psi \in \mathcal{S}(G)$ is called LP-admissible if for all $g \in \mathcal{Z}(G)$,
\begin{equation} \label{eqn:LP_ZG}
 g = \lim_{N \to \infty} \sum_{|j| \le N} g \ast \psi_j^* \ast \psi_j
\end{equation} holds, with convergence in the Schwartz space topology.
 Duality entails the convergence
\begin{equation} \label{eqn:LP_dist}
 u = \lim_{N \to \infty} \sum_{|j| \le N} u \ast \psi_j^* \ast \psi_j~.
\end{equation} for all $u \in \mathcal{S}'(G)/\mathcal{P}$.
\end{definition}

The following lemma yields the chief construction of LP-admissible functions
\begin{lemma} \label{lem:construct_LP_adm}
Let $\hat\phi$ be a function in $C^\infty$ with support in $[0,4]$
such that $0\leq \hat\phi\leq 1$ and  $\phi\equiv 1$
  on $[0,1/4]$. Take
  $\hat \psi(\xi)= \sqrt{\hat\phi(2^{-2}\xi) -\hat\phi(\xi)}$. Thus, $\hat\psi \in C_c^\infty(\RR^+)$,
with support in the  interval $[1/4, 4]$, and
\begin{equation} \label{eqn:LP_condition_Fourier}
 1= \sum_{j\in \ZZ} |\hat\psi(2^{2j}\xi)|^2\quad a.e ~~
\end{equation}
Pick a sub-Laplacian $L$, and let $\psi$ denote the convolution kernel associated to the bounded left-invariant operator $\widehat{\psi}(L)$. Then $\psi$  is LP-admissible, with $\psi \in \mathcal{Z}(G)$.
\end{lemma}

\begin{proof}
Let us first comment on the properties of $\psi$ that are immediate from the construction via spectral calculus: $\psi \in \mathcal{S}(G)$ follows from \cite{Hu}, and vanishing moments by \cite[Proposition 1]{gm1}.

Now let $g\in \mathcal{Z}(G)$. First note that  $2$-homogeneity of $L$ implies that the convolution kernel associated to $\widehat{\psi}(2^{-2j} \cdot )(L)$ coincides with $\psi_j$. Then, by the spectral theorem and (\ref{eqn:LP_condition_Fourier}),
\begin{align}\label{conver-strong}
g = \sum_{j\in \ZZ}   [\hat\psi(2^{-2j} \cdot )(L)]^\ast  \circ
[\hat\psi(2^{-2j} \cdot)(L)] g   =  \sum_{j \in \ZZ} g \ast \psi_j^* \ast \psi_j \end{align}
 holds in $L^2$-norm.

For any positive integer $N$
\begin{align}\notag
 \sum_{\mid j\mid\leq N} g\ast\psi_{j}^\ast\ast\psi_j  = g\ast D_{2^{N+1}} \phi - g\ast D_{2^{-N}} \phi ~,
 \end{align} where $\phi \in \mathcal{S}(G)$ is the convolution kernel of $\widehat{\phi}(L)$.
Since $\phi$ is a Schwartz function, it follows by \cite[Proposition (1.49)]{FollandStein82} that $g \ast D_{2^{N+1}} \phi  \to c_\phi g$, for $N \to \infty$, for all $g \in \mathcal{S}(G)$, with convergence in $\mathcal{S}(G)$ and a suitable constant $c_\phi$.

We next show that $g \ast D_t f \to 0$ in $\mathcal{S}(G)$, as $t \to 0$, for any $f \in \mathcal{S}(G)$. Fix a multi-index $I$ and $N,k  \in \mathbb{N}$ with $k \ge N$. Then left-invariance and homogeneity of $Y^I$ yield
\begin{eqnarray*}
 | Y^I (g \ast D_t f)(x)| & = & t^{d(I)} |g \ast D_t (Y^I f)(x)|  \\
 & \le & C_{f,g} t^{k+Q+d(I)}(1+|tx|)^{-N} \\
& \le & C_{f,g} t^{k+Q+d(I)-N} (1+|x|)^{-N}~.
\end{eqnarray*} Here the first inequality is an application of (\ref{eqn:decay_conv1});  the constant $C_{f,g}$ can be estimated in terms of $|f|_M,|g|_M$, for $M$ sufficiently large. But this proves $g \ast D_t f \to 0$ in the Schwartz topology.

Summarizing, $\sum_{|j| \le N} g\ast\psi_{j}^\ast\ast\psi_j \to c_\phi g$ in $\mathcal{S}(G)$, and in addition by
(\ref{conver-strong}), $\sum_{|j| \le N} g\ast\psi_{j}^\ast\ast\psi_j \to g$ in ${\rm L}^2$, whence $c_\phi=1$ follows.
\end{proof}

Note that an LP-admissible function $\psi$ as constructed in \ref{lem:construct_LP_adm} fulfills the convenient equality
\begin{equation} \label{eqn:conv_vanish}
\forall j, l \in \mathbb{Z}~:~|j-l| > 1 \Rightarrow \psi_j^* \ast \psi_l =0 ~,
\end{equation} which follows from $[\hat\psi(2^{-2j} \cdot)(L)]\circ [\hat\psi(2^{-2l} \cdot)(L)]= 0$.

\begin{remark} \label{rem:dense_image}
 By spectral calculus, we find that $\psi = L^k g_k$, with $g_k \in \mathcal Z(G)$. 
 In particular, the decomposition
\begin{eqnarray*}
 f & = & \lim_{N \to \infty} \sum_{|j| \le N} f \ast \psi_j^* \ast D_{2^j} L^k(g_k)  \\
 & = & \lim_{N \to \infty} L^k \left( \sum_{|j| \le N} f \ast \psi_j^* \ast 2^{-kj} D_{2^j} g_k \right)
\end{eqnarray*}
shows that $L^k(\mathcal{Z}(G)) \subset \mathcal{Z}(G)$ is dense.
\end{remark}

 We now associate a scale of homogeneous Besov  spaces to the function $\psi$.
\begin{definition} \label{defn:B_psi}
 Let $\psi \in \mathcal{Z}(G)$ be LP-admissible, let $1 \le p \le \infty$, $1\le q \le  \infty$, and $s \in \mathbb{R}$. The {\bf homogeneous Besov space associated to $\psi$} is defined as
\begin{equation}
 \label{eqn:def_B_psi}
\dot{B}_{p,q}^{s,\psi} = \left\{ u \in \mathcal{S}'(G)/\mathcal{P} : \left\{ 2^{js} \| u \ast \psi_j^* \|_p \right\}_{j \in \mathbb{Z}} \in \ell^q(\mathbb{Z}) \right\}~,
\end{equation}
with associated norm
\[
 \| u \|_{\dot{B}_{p,q}^{s,\psi}} = \left\|  \left\{ 2^{js} \| u \ast \psi_j^* \|_p \right\}_{j \in \mathbb{Z}}  \right\|_{\ell^q(\mathbb{Z})}~.
\]
\end{definition}

\begin{remark} The definition relies on the conventions regarding ${\rm L}^p$-norms of distributions (modulo polynomials), as outlined in Remark \ref{rem:LP_norm_dist}.
Definiteness of the Besov norm holds because of (\ref{eqn:LP_dist}).
\end{remark}

The combination of Lemma \ref{lem:construct_LP_adm} with Definition \ref{defn:B_psi} shows that we cover the homogeneous Besov spaces defined in the usual manner via the spectral calculus of sub-Laplacians.
Hence the following theorem implies in particular that different sub-Laplacians yield the same homogeneous Besov spaces (at least within the range of sub-Laplacians that we consider).

\begin{theorem}
 \label{thm:Besov_norm_equiv}
 Let $\psi^1,\psi^2 \in \mathcal{Z}(G)$ be LP-admissible. Let $s \in \mathbb{R}$ and $1 \le p,q \le \infty$.  Then  $\dot{B}_{p,q}^{s,\psi^1} = \dot{B}_{p,q}^{s,\psi^2}$, with equivalent norms.
\end{theorem}
\begin{proof}
 It is sufficient to prove the norm equivalence, and here symmetry with respect to $\psi^1$ and $\psi^2$ immediately reduces the proof to showing, for a suitable constant $C>0$,
\begin{equation} \label{eqn:half_norm_equiv}
  \forall u \in \mathcal{S}'(G)/\mathcal{P} : \| u \|_{\dot{B}_{p,q}^{s,\psi^1}} \le C \| u \|_{\dot{B}_{p,q}^{s,\psi^2}} ~,
\end{equation} in the extended sense that the left-hand side is finite whenever the right-hand side is. Hence assume that $u \in \dot{B}_{p,q}^{s,\psi^2}$; otherwise, there is nothing to show.  In the following, let $\psi_{i,j} = D_{2^j} \psi^i$ (i=1,2).

By LP-admissibility of $\psi^2$,
\[
 u = \lim_{N \to \infty} \sum_{|j| \le N} u \ast \psi_{2,j}^* \ast \psi_{2,j}~,
\] with convergence in $\mathcal{S}'(G)/\mathcal{P}$.
Accordingly,
\begin{equation} \label{eqn:u_ast_eta}
 u \ast \psi_{1,\ell}^*= \lim_{N \to \infty} \sum_{|j| \le N} u \ast  \psi_{2,j}^* \ast  \psi_{2,j} \ast \psi_{1,\ell}^* ~,
\end{equation} where the convergence on the right-hand side holds in $\mathcal{S}'(G)$, by \ref{lem:conv_SmP}. We next show that the right-hand side also converges in ${\rm L}^p$. For this purpose, we observe that
\[
 \|  \psi_{2,j} \ast \psi_{1,\ell}^* \|_1 = \| D_{2^j} (\psi^2 \ast D_{2^{\ell-j}} \psi_1^{1*}) \|_1 = \| \psi^2 \ast D_{2^{\ell-j}} \psi^{1*} \|_1
 \le C 2^{-|\ell- j|k}~,
\] where $k>s$ is a fixed integer. For $\ell-j \ge 0$, this follows directly from (\ref{eqn:decay_conv2_norm}),  using $\psi^1,\psi^2 \in \mathcal{S}(G)$,  and vanishing moments of $\psi^1$, whereas for $\ell-j <0$, the vanishing moments of $\psi^2$ allow to apply (\ref{eqn:decay_conv1_norm}).

Using Young's inequality, we estimate with $C$ from above that
\begin{eqnarray}
\nonumber \sum_{j \in \mathbb{Z}}  \| u \ast  \psi_{2,j}^* \ast \psi_{2,j} \ast \psi_{1,\ell}^* \|_p & \le &  \sum_{j \in \mathbb{Z}} \| u \ast \psi_{2,j}^*\|_p \| \psi_j \ast \psi_\ell^* \|_1 \\ & \le & C  \| u \ast \psi_{2,j}^*\|_p  2^{-|j-\ell|k} \\
 & \le &  \label{eqn:est_sumj} C \sum_{j  \in \mathbb{Z}} 2^{js} \| u \ast \psi_{2,j}^*\|_p   2^{-|j-\ell|k-js}~.
\end{eqnarray}
Next observe that
\begin{equation} \label{eqn:matrix_decay} 2^{-|j-\ell|k - js}  = 2^{-\ell s} \cdot \left\{ \begin{array}{cc} 2^{-|j-\ell|(k+s) } & j \ge \ell \\ 2^{-|j-\ell| (k-s)} & j < \ell \end{array} \right. \le 2^{-\ell s} 2^{-|j-\ell|(k-|s|)}~.
\end{equation}
By assumption, the sequence $(2^{js} \| u \ast \psi_{j,2}^* \|_p)_{j \in \mathbb{Z}}$ is in $\ell^q$, in particular bounded. Therefore, $k-|s|>0$ yields that (\ref{eqn:est_sumj}) converges. But then the right-hand side of (\ref{eqn:u_ast_eta}) converges unconditionally with respect to $\| \cdot \|_p$. This limit coincides with the $\mathcal{S}'(G)/\mathcal{P}$-limit $u \ast \psi_{1,\ell}^*$ (which because of $\psi_{1,\ell}^* \in \mathcal{Z}(G)$ is even a $\mathcal{S}'(G)$-limit), yielding $u \ast \psi_{1,\ell}^* \in {\rm L}^p(G)$, with
\begin{eqnarray*}
 2^{\ell s} \| u \ast \psi_{1,\ell}^* \|_p & \le &  2^{\ell s} \sum_{j \in \mathbb{Z}}  \| u \ast  \psi_{2,j}^* \ast \psi_{2,j} \ast \psi_\ell^* \|_p \\
 & \le & C_3 2^{\ell s} \sum_{j \in \mathbb{Z}} 2^{js} \| u \ast \psi_{2,j}^* \|_p 2^{-|j-\ell| (k-|s|)}~.
\end{eqnarray*}
Now an application of Young's inequality for convolution over $\mathbb{Z}$, again using $k - |s|>0$, provides (\ref{eqn:half_norm_equiv}).
\end{proof}

As a consequence, we write $\dot{B}_{p,q}^s = \dot{B}_{p,q}^{s,\psi}$, for any LP-admissible $\psi \in \mathcal{Z}(G)$. These spaces coincide with the homogeneous Besov spaces for the Heisenberg group in \cite{Ba}, and with the usual definitions in the case $G=\mathbb{R}^n$.

In the remainder of the section we note some functional-analytic properties of Besov spaces and Littlewood-Paley-decompositions for later use.
\begin{lemma} \label{lem:cont_embed}
 For all $1 \le p,q \le \infty$ and all $s \in \mathbb{R}$, one has continuous inclusion maps $\mathcal{Z}(G) \hookrightarrow {\dot B}_{p,q}^s \hookrightarrow S'(G)/\mathcal{P}$, as well as $\mathcal{Z}(G) \hookrightarrow {\dot B}_{p,q}^{s*}$, where the latter denotes the dual of ${\dot B}_{p,q}^s$. For $p,q < \infty$, $\mathcal{Z}(G) \subset {\dot B}_{p,q}^s$ is dense.
\end{lemma}
\begin{prf} We pick $\psi$ as in Lemma \ref{lem:construct_LP_adm} and define $\Delta_j g= g\ast \psi_j^*$ for $g \in \mathcal{S}'(G)$. 
For the inclusion $\mathcal{Z}(G) \subset \dot{B}_{p,q}^s$, note
that (\ref{eqn:decay_conv1_norm}) and (\ref{eqn:decay_conv2_norm})
allow to estimate for all $g \in \mathcal{Z}(G)$ and $k \in
\mathbb{N}$ that
\[
 \| \Delta_j g \|_p \le C_{k} 2^{-|j|k}
\] Here the constant $C_k$ is a suitable multiple of $|g|_M$, for $M=M(k)$ sufficiently large. But this implies that $\mathcal{Z}(G) \subset \dot{B}_{p,q}^s$ continuously.

For the other embedding, repeated applications of H\"older's inequality yield the estimate
\begin{eqnarray*}
 |\langle f,g \rangle| & = & |\sum_{j \in \mathbb{Z}} \langle f, g \ast \psi_j^* \ast \psi\rangle| \\
 &  \le &  \sum
_{j \in \mathbb{Z}} | \langle f \ast \psi_j^*, g \ast \psi_j^* \rangle |  \\
 & \le &  \sum_{j \in \mathbb{Z}} \| f \ast \psi_j^* \|_{p'} ~  \| g \ast \psi_j^* \|_{p} \\
 & = & \sum_{j \in \mathbb{Z}} (2^{-js} \| f \ast \psi_j^*\|_{p'})  (2^{js} \| f \ast \psi_j^* \|_p )\\
 & \le & \| f \|_{p',q'}^{-s}~ \| g \|_{p,q}^{s}
\end{eqnarray*} valid for all $f \in \mathcal{Z}(G) \subset {\dot B}_{p',q'}^{-s}$ and $g \in {\dot B}_{p,q}^s$. Here $p',q'$ are the conjugate exponents of $p,q$, respectively. But this estimate implies continuity of the embeddings ${\dot B}_{p,q}^s \subset \mathcal{S}'(G)/\mathcal{P}$ and $\mathcal{Z}(G) \subset {\dot B}_{p,q}^{s*}$.

For the density statement, let $u \in {\dot B}_{p,q}^s$, and $\epsilon>0$. For convenience, we pick $\psi$ according to Lemma \ref{lem:construct_LP_adm}. Since $q< \infty$, there exists $N \in \mathbb{N}$ such that
\[
 \sum_{|j| >N-1} 2^{jsq} \| \Delta_j u \|_p^q < \epsilon~.
\] Next define
\[ K_N = \sum_{|j| \le N} \psi_j^* \ast \psi_j  = D_{2^{N+1}} \phi - D_{2^{-N}} \phi ~.\]
Let $w = u \ast K_N$. By assumption on $u$ and Young's inequality, $w \in {\rm L}^p(G)$, and since $p< \infty$, there exists $g \in \mathcal{S}(G)$ with $\| w - g \|_p < \epsilon^{1/q}$. Let $f = g \ast K_N$, then $f \in \mathcal{Z}(G)$, and for $j \in \mathbb{Z}$,
\begin{eqnarray*}
 \| \Delta_j (u-f) \|_p  &  = &  \| (u - f) \ast \psi_j^* \|_p  \\
 & \le &  \| u \ast \psi_j^* - u \ast K_N \ast \psi_j^* \|_p + \| w \ast \psi_j^* - g \ast K_N \ast \psi_j^* \|_p
\end{eqnarray*}
For $|j| \le  N -1 $, the construction of $\psi_j$ and $K_N$ implies  that  $K_N \ast \psi_j^* = \psi_j^*$, whereas for $|j| > N+1$, one has $K_N \ast \psi_j^* = 0$. As a consequence, one finds for $|j| < N-1$
\[
 \| \Delta_j (u-f) \|_p \le \| w - g \|_p \| \psi_j^* \|_1 = \| w -g \|_p \| \psi \|_1 <  \epsilon^{1/q} \| \psi \|_1~,
\] and for $|j| > N+1$
\[
 \| \Delta_j (u-f) \|_p \le \| u \ast \psi_j^* \|_p < \epsilon^{1/q}~.
\]
For $|~|j|~-N| \le 1$, one finds
\[
 \| \Delta_j (u-f) \|_p \le C\epsilon^{1/q}
\] with some constant $C>0$ depending only on $\psi$.
For instance, for $j=N$,
\begin{eqnarray*}  \| \Delta_j (u-f) \|_p & \le & \| u \ast \psi_N^* - u \ast (\psi_{N-1}^* \ast \psi_{N-1} + \psi_N^ * \ast \psi_N ) \ast \psi_N^* \|_p  \\ &  + & \| w  \ast \psi_N^* - g  \ast  (\psi_{N-1}^* \ast \psi_{N-1} + \psi_N^ * \ast\psi_N) \ast \psi_N^*  \|_p~.
\end{eqnarray*}
A straigthforward application of triangle and Young's inequality yields:
\[
 \| u \ast \psi_N^* - u \ast (\psi_{N-1}^* \ast \psi_{N-1} + \psi_N^ * \ast \psi_N ) \ast \psi_N^* \|_p   \le \| u \ast\psi_N^* \|_p (1+2 \| \psi^* \ast \psi \|_1) < \epsilon^{1/q}(1+2 \| \psi^* \ast \psi \|_1) ~.
\] Similar considerations applied to $w = u \ast K_N$ yield
\begin{eqnarray*}
\lefteqn{  \| w  \ast \psi_N^* - g  \ast  (\psi_{N-1}^* \ast \psi_{N-1} + \psi_N^ * \ast\psi_N) \ast \psi_N^*  \|_p} \\  & \le &
2 \| u \ast \psi_N^* \|_p \| \psi^* \ast \psi \|_1 + 2 \| g \ast \psi_N^* \|_p \| \psi^* \ast \psi \|_1  \\
& \le & 2 \epsilon^{1/q} \| \psi^* \ast \psi \|_1 + 2 (\| w \ast \psi_N^* \|_p + \| (w-g) \ast \psi_N^* \|_p)   \| \psi^* \ast \psi \|_1 \\
& \le & (4 \| \psi^* \ast \psi \|_1 + \| \psi^* \ast \psi \|_1 \|\psi \|_1) \epsilon^{1/q}.
\end{eqnarray*}
Now summation over $j$ yields
\[
 \| u-f \|_{{\dot B}_{p,q}^s} \le C' \epsilon~,
\] as desired.
\end{prf}

\begin{remark} \label{rem:LP_finite}
 Let $\psi$ as in Lemma \ref{lem:construct_LP_adm}. As a byproduct of the proof, we note that the space
\[
 \mathcal{D} = \{ f \ast K_N ~:~f \in \mathcal{S}(G), N \in \mathbb{N} \}
\] is dense in $\mathcal{Z}(G)$ as well as ${\dot B}_{p,q}^s$, if $p,q < \infty$. In $\mathcal{D}$, the decomposition
\[
 g = \sum_{j \in \mathbb{Z}} g \ast \psi_j^* \ast \psi_j
\] holds with finitely many nonzero terms.
\end{remark}

We next extend the Littlewood-Paley decomposition to the elements of the Besov space. For simplicity, we prove the result only for certain LP-admissible functions.
  \begin{proposition} \label{prop:decomp_conv_Besov} Let $1 \le p,q < \infty$, and $\psi \in \mathcal{Z}(G)$ an LP-admissible vector constructed via Lemma \ref{lem:construct_LP_adm}.
  Then the decomposition (\ref{eqn:LP_ZG}) converges for all $g \in {\dot B}_{p,q}^s$ in the Besov space norm.
   \end{proposition}
 \begin{proof}
Consider the operators $\Sigma_N : {\dot B}_{p,q}^s \to {\dot B}_{p,q}^s$,
\[
 \Sigma_N u = \sum_{|j| \le N} u \ast \psi_j^* \ast \psi_j~.
\] By suitably adapting the arguments proving the density statement of Lemma \ref{lem:cont_embed}, it is easy to see that the family of operators $(\Sigma_N)_{N \in \mathbb{N}}$ is bounded in the operator norm.
 As noted in \ref{rem:LP_finite}, the $\Sigma_N$ strongly converges to the identity operator on a dense subspace. But then boundedness of the family implies strong convergence everywhere.
\end{proof}

A further class of spaces for which the decomposition converges is ${\rm L}^p$:
\begin{proposition} \label{prop:decomp_conv_Lp}
Let $1 < p < \infty$, and $\psi \in \mathcal{Z}(G)$ an LP-admissible vector constructed via Lemma \ref{lem:construct_LP_adm}. Then the decomposition (\ref{eqn:LP_ZG}) converges with respect to $\| \cdot \|_p$, for all $g \in {\rm L}^p(G)$.
\end{proposition}
\begin{proof}
 Let the operator family $(\Sigma_N)_{N \in \mathbb{N}}$ be defined as in the previous proof. Then $\Sigma_N f = g\ast D_{2^{N+1}} \phi - g\ast D_{2^{-N}} \phi $, and Young's inequality implies that the sequence of operators is norm-bounded. It therefore suffices to prove the desired convergence on the dense subspace $\mathcal{S}(G)$. By \cite[1.20]{FollandStein82},
 $g \ast D_{2^{N+1}} \phi \to c_\phi g$.
Furthermore, for $N \in \mathbb{N}$,
\begin{eqnarray*}
 \left( g \ast D_{2^{-N}} \phi\right) (x) & = & 2^{-NQ} \int_G g(y) \phi(2^{-N}(y^{-1}x)) dy \\
 & = & \int_G g(2^{N} y) \phi(y^{-1} \cdot 2^{-N} x) dy \\
& = & 2^{-NQ} \left( D_{2^N} g \ast \phi \right) (2^{-N} x) ~,
\end{eqnarray*}
and thus
\begin{eqnarray*}
 \| g \ast D_{2^{-NQ}} \phi \|_p & = & 2^{-NQ} \left( \int_{G} \left| (D_{2^N} g \ast \phi) (2^{-N} x)\right|^p  dx \right)^{1/p} \\ & = & 2^{-NQ+NQ/p}   \| D_{2^N} g \ast \phi\|_p~.
\end{eqnarray*}
Again by \cite[1.20]{FollandStein82}, $(D_{2^N} g \ast \phi) \to c_g \phi$, in particular,
\[
  2^{-NQ+NQ/p}   \| D_{2^N} g \ast \phi\|_p \to 0 \mbox{ as } N \to \infty.
\] Hence $\Sigma_N g \to c_\phi g$, and the case $p=2$ yields $c_\phi=1$.
\end{proof}

 \begin{theorem}  ${\dot B}_{p ,q}^{s}$ is a Banach space.
 \end{theorem}
 \begin{proof}
Completeness is the only issue here. Again, we pick $\psi \in \mathcal{Z}(G)$ an LP-admissible vector via Lemma \ref{lem:construct_LP_adm}. Suppose that $\{ u_n \}_{n \in \mathbb{N}} \subset \dot{B}_{p,q}^s$ is a Cauchy sequence.
As a consequence, one has in particular, for all $j \in \mathbb{Z}$, that $ \{ u_n \ast \psi_j^* \}_{n \in \mathbb{N}} \subset {\rm L}^p(G)$ is a Cauchy sequence, hence $u_n \ast \psi_j^* \to v_j$, for a suitable $v_j \in {\rm L}^p(G)$. Furthermore, the Cauchy property of $\{ u_n \}_{n \in \mathbb{N}} \subset  \dot{B}_{p,q}^s$ implies that
\[
 \left\{ \left\{ 2^{js} \| u_n \ast \psi_{j}^* \|_p \right\}_{j \in \mathbb{Z}} \right\}_{n \in \mathbb{N}} \subset \ell^q(\mathbb{Z})
\] is a Cauchy sequence. On the other hand, the sequence converges pointwise to $\left\{ 2^{js} \| v_j \|_p \right\}_j$, whence
\begin{equation} \label{eqn:limit_lq}
\sum_{j \in \mathbb{Z}}  2^{jsq} \| v_j \|_p^q < \infty~.
\end{equation}
We define
\[
 u = \lim_{M \to \infty} \sum_{|j| \le M}  v_j \ast \psi_j~.
\] Now,  using  (\ref{eqn:limit_lq}) and $\mathcal{Z}(G) \subset {\dot B}_{p',q'}^{-s}$, where $p',q'$ are the conjugate exponents of $p,q$, respectively, a straightforward calculation as in the proof of Lemma \ref{lem:cont_embed} shows that the sum defining $u$ converges in $\mathcal{S}'(G)/\mathcal{P}$. Furthermore, (\ref{eqn:limit_lq}) and (\ref{eqn:conv_vanish}) easily imply that $u \in \dot{B}_{p,q}^s$. Finally, for the proof of $u_n \to u$, we employ  (\ref{eqn:conv_vanish}) together with the equality $\psi_j^* = \sum_{|l-j| \le 1} \psi_l^* \ast \psi_l \ast \psi_j^*$, to show that
\begin{eqnarray*}
 \| (u_n - u) \ast \psi_j \|_p & = &  \| u_n \ast \psi_j - \sum_{|l-j| \le 1}  v_l \ast \psi_l \ast \psi_j^* \|_p \\
 & \le & \sum_{|l-j|\le 1} \| (u_n \ast \psi_l^*- v_l) \ast \psi_l \ast \psi_j^* \|_p \\
 & \le & \sum_{|l-j| \le 1} \| u_n \ast \psi_l^* - v_l \|_p \| \psi_l \ast \psi_j^* \|_1 \to 0~,\mbox{ as } n \to \infty~.
\end{eqnarray*}
 Summarizing, the sequence $\left\{  \left\{ 2^{js} \| (u_n - u) \ast \psi_j^* \|_p \right\}_{j \in \mathbb{Z}} \right\}_{n \in \mathbb{N}} \in \ell^q(\mathbb{N})$ is a Cauchy sequence, converging pointwise to $0$.  But then $\| u_n - u\|_{\dot{B}_{p,q}^s} \to 0$ follows.
  \end{proof}

\section{Characterization via Continuous Wavelet Transform}

The following definition can be viewed as a continuous-scale analog of LP-admissibility.

\begin{definition}
 $\psi \in \mathcal{S}(G)$ is called {\bf $\mathcal{Z}$-admissible}, if for all $f \in \mathcal{Z}(G)$,
\[
  f = \lim_{\epsilon \to 0, A \to \infty} \int_\epsilon^A f \ast D_a(\psi^* \ast \psi) \frac{da}{a}~
\] holds with convergence in the Schwartz topology.
\end{definition}

The next theorem reveals a large class of $\mathcal{Z}$-admissible
wavelets. In fact, all the wavelets studied in \cite{gm1} are also
$\mathcal{Z}$-admissible in the sense considered here. Its proof is
an adaptation of the argument showing \cite[Theorem 1]{gm1}.
\begin{theorem} \label{thm:adm_cond_CWT}
 Let $\widehat{h} \in \mathcal{S}(\mathbb{R}^+)$, and let $\psi $ be the distribution kernel associated to the operator $L \widehat{h}(L)$.  Then $\psi$ is $\mathcal{Z}$-admissible up to normalization.
\end{theorem}

\begin{proof}
The main idea of the proof is to write, for $f \in \mathcal{Z}(G)$,
\begin{eqnarray*} \int_\epsilon^A f \ast D_a(\psi^* \ast \psi) \frac{da}{a} & = & f \ast \int_{\epsilon}^A D_a(\psi^* \ast \psi) \frac{da}{a} \\ & = & f \ast D_A g - f \ast D_\epsilon g~,\end{eqnarray*}
with suitable $g \in \mathcal{S}(G)$. Once this is established, $f
\ast D_A g \to c_g f$ for $A \to \infty$ follows by
\cite[Proposition (1.49)]{FollandStein82}, with convergence in the
Schwartz topology. Moreover, $f \in \mathcal{Z}(G)$ entails that $f
\ast D_\epsilon g \to 0$ in the Schwartz topology: Given any $N>0$
and $I \in \mathbb{N}^n$ with associated left-invariant differential
operator $Y^I$, we can employ (\ref{eqn:decay_conv1})  to estimate
\begin{eqnarray*}
 \sup_{x \in G} (1+|x|)^N  \left| (Y^I f \ast D_{\epsilon} g)(x) \right| & = & \sup_{x \in G} (1+|x|)^N \epsilon{Q+d(I)}
\left| f \ast D_{\epsilon} (Y^I g)(x) \right| \\
 & \le & C \sup_{x \in G} (1+|x|)^N \epsilon^{Q+d(I)+k}  (1+|\epsilon x|)^{-M} \\
 & \le & C \sup_{x \in G} (1+|x|)^{N-M} \epsilon^{Q+d(I)+k+M} ~,
\end{eqnarray*} which converges to zero for $\epsilon \to 0$, as soon as $M\ge N$ and $k > M-Q-d(I)$. But this implies $f \ast D_\epsilon g \to 0$ in $\mathcal{S}(G)$, by \cite{FollandStein82}.

Thus it remains to construct $g$. To this end, define
\begin{equation}
 \widehat{g}(\xi) = -\frac{1}{2} \int_{\xi}^\infty  a|\widehat{h}(a^2)|^2 da~,
\end{equation} which is clearly in $\mathcal{S}(\mathbb{R}^+)$, and let $g$ denote the associated convolution kernel of $ \widehat{g}(L)$. By the definition, $g\in \mathcal{S}(G)$. Let $\varphi_1,\varphi_2$ be in $\mathcal{S}(G)$, and let $d\lambda_{\varphi_1,\varphi_2}$ denote the scalar-valued Borel measure associated to $\varphi_1,\varphi_2$ by the spectral measure. Then, by spectral calculus and the invariance properties of $da/a$,
\begin{eqnarray*}
 \langle \int_{\epsilon}^A \varphi_1 \ast D_a (\psi^* \ast \psi) f \frac{da}{a}, \varphi_2 \rangle & = &  \int_{0}^\infty \int_\epsilon^A (a^2 \xi)^2 |\hat h(a^2 \xi)|^2 \frac{da}{a} d\lambda_{\varphi_1,\varphi_2}(\xi)  \\
& = & \frac{1}{2}\int_{0}^\infty \int_{\epsilon^2 \xi}^{A^2 \xi} a |\hat h(a^2 \xi)|^2 da d\lambda_{\varphi_1,\varphi_2}(\xi) \\
& = & \int_{0}^\infty \widehat{g}(A^2 \xi) - \widehat{g}(\epsilon^2 \xi)  d\lambda_{\varphi_1,\varphi_2}(\xi) \\
& = & \langle \varphi_1 \ast (D_A g - D_\epsilon g), \varphi_2 \rangle ~,
\end{eqnarray*} as desired.
\end{proof}

Hence, by \cite[Corollary 1]{gm1}:
\begin{corollary}
 \begin{enumerate}
  \item [(a)] There exist $\mathcal{Z}$-admissible $\psi \in \mathcal{Z}(G)$.
  \item[(b)] There exist $\mathcal{Z}$-admissible $\psi \in C_c^\infty(G)$ with vanishing moments of arbitrary finite order.
 \end{enumerate}
\end{corollary}

Given a tempered distribution $u \in \mathcal{S}'(G)/\mathcal{P}$
and a $\mathcal{Z}(G)$-admissible function $\psi$, the continuous
wavelet transform of $u$ is the family $(u \ast D_a \psi^*)_{a > 0}$
of convolution products. We will now prove a characterization of
Besov spaces in terms of the continuous wavelet transform. 

Another popular candidate for defining scales of Besov spaces is the
heat semigroup; see e.g. \cite{Saka79} for the inhomogeneous case on
stratified groups, or   rather  \cite{BuBe} for the general treatment.
In our setting, the heat semigroup associated to the sub-Laplacian
is given by right convolution with $h_t(x) = D_t h(x)$, where $h$ is the kernel of $\widehat{h}(L)$ with $\widehat{h}(\xi) = e^{-\xi}$. Theorem
\ref{thm:adm_cond_CWT} implies that $\psi = L^k h$ is
$\mathcal{Z}$-admissible; it can be viewed as an analog of the
well-known {\bf Mexican Hat} wavelet. (This wavelet  on the general groups  was studied for the first time in \cite{gm1}.) 
The wavelet transform of $f
\in \mathcal{S}'(G)$ associated to $\psi$ is then very closely
related to the $k$-fold time derivative of the solution to the heat
equation with initial condition $f$: By choice of $h$,
\[
u(x,t) = (f \ast D_t h)(x)~,
\] denotes the solution of the heat equation associated to $L$, with initial condition $f$. A formal calculation using left invariance of $L$ then yields
\[
\partial_t^k u = L^k (f \ast D_t h) = f \ast L^k (D_t h) = t^{2k} f
\ast D_t \psi^* ~.
\]
Thus the following theorem also implies a characterization of Besov
spaces in terms of the heat semigroup.

\begin{theorem} \label{thm:Char_CWT}
Let $\psi \in \mathcal{S}(G)$ be $\mathcal{Z}$-admissible, with vanishing moments of order $k$. Then, for all $s \in \mathbb{R}$ with $|s|< k$, and all $1 \le p < \infty$, $1 \le q \le \infty$, the following norm equivalence holds:
\begin{equation} \label{eqn:equiv_CWT_Besov}
\forall u \in \mathcal{S}'(G)/\mathcal{P}~:~ \| u \|_{\dot{B}_{p,q}^s} \asymp \left\| a \mapsto  a^{s} \| u \ast D_a \psi^* \|_p \right\|_{{\rm L}^q(\mathbb{R}^+; \frac{da}{a})}  \end{equation}
 Here the norm equivalence is understood in the extended sense that one side is finite iff the other side is. If $\psi \in \mathcal{Z}(G)$, the equivalence is also valid for the case $p=\infty$.
\end{theorem}
\begin{proof}
The strategy consists in adapting the proof of Theorem \ref{thm:Besov_norm_equiv} to the setting where one summation over scales is replaced by integration. This time however, we have to deal with both directions of the norm equivalence.
In the following estimates, the symbol $C$ denotes a constant that may change from line to line, but in a way that is independent of $u \in S'(G)$.

Let us first assume that
\[
 \int_{\mathbb{R}} a^{sq} \| u \ast D_a \psi^* \|_p^q \frac{da}{a} < \infty~,
\] for $u \in \mathcal{S}'(G)/\mathcal{P}$, $1 \le p,q \le \infty$, for a $\mathcal{Z}$-admissible function $\psi \in S(G)$ with $k_\psi>|s|$ vanishing moments ($\psi \in \mathcal{Z}(G)$, if $p=\infty$). Let $\varphi \in \mathcal{Z}(G)$ be LP-admissible. Then, for all $j \in \mathbb{Z}$,
\begin{equation} \label{eqn:decomp_u_varphi}
 u \ast \varphi_j ^*= \lim_{\epsilon \to 0, A \to \infty} \int_{\epsilon}^A u \ast D_a \psi^* \ast D_a \psi \ast \varphi_j^* \frac{da}{a}
\end{equation} holds in $\mathcal{S}'(G)$, by \ref{lem:conv_SmP}.

We next prove that the right-hand side of (\ref{eqn:decomp_u_varphi}) converges in ${\rm L}^p$. For this purpose, introduce
\[
 c_j = \int_{0}^\infty \| u \ast D_a \psi^* \ast D_a \psi \ast \varphi_j^* \|_p \frac{da}{a}~.
\]

 We estimate
\begin{eqnarray}
 \nonumber
 c_j & \le & \int_{0}^\infty \| u \ast D_a \psi^* \|_p \| D_a \psi \ast \varphi_j^* \|_1 \frac{da}{a} \\
 & = & \label{eqn:anschluss_infty} \int_{1}^2 \sum_{\ell \in \mathbb{Z}} \| u \ast D_{a 2^\ell} \psi^* \|_p \| D_{a 2^\ell} \psi \ast \varphi_j^* \|_1 \frac{da}{a} \\
& \le & \label{eqn:anschluss_q} \left( \int_{1}^2 \left( \sum_{\ell \in \mathbb{Z}} \| u \ast D_{a 2^\ell} \psi^* \|_p \| D_{a 2^\ell} \psi \ast \varphi_j^* \|_1 \right)^q \frac{da}{a} \right)^{1/q} \log(2)^{1/q'}~,
\end{eqnarray} where we used that $da/a$ is scaling-invariant. Note that the last inequality is H\"older's inequality for $q< \infty$. In this case, taking $q$th powers and summing over $j$ yields
\begin{equation} \label{eqn:onedir_integrand}
 \sum_{j \in \mathbb{Z}} 2^{jsq} c_j^q \le C
\int_1^2 \sum_{j \in \mathbb{Z}} 2^{jsq} \left(  \sum_{\ell \in \mathbb{Z}} \| u \ast D_{a 2^\ell} \psi^* \|_p \| D_{a 2^\ell} \psi \ast \varphi_j^* \|_1 \right)^q \frac{da}{a}~.
\end{equation}

 Using vanishing moments and Schwartz properties of $\psi$ and $\varphi$, we can now employ  (\ref{eqn:decay_conv1_norm}) and (\ref{eqn:decay_conv2_norm}) to obtain
\begin{equation} \label{eqn:l1_kernel}
 \| D_{a 2^\ell} \psi \ast \varphi_j^* \|_1 \le C 2^{-|j-\ell|k}~,
\end{equation} with a constant independent of $a \in [1,2]$. But then, since $k>|s|$, we may proceed just as in the proof of \ref{thm:Besov_norm_equiv} to estimate the integrand in (\ref{eqn:onedir_integrand}) via
\begin{equation}
 \sum_{j \in \mathbb{Z}} 2^{jsq} \left(  \sum_{\ell \in \mathbb{Z}} \| u \ast D_{a 2^\ell} \psi^* \|_p \| D_{a 2^\ell} \psi \ast \varphi_j^* \|_1 \right)^q \le C\sum_{\ell \in \mathbb{Z}} 2^{\ell s q} \| u \ast D_{a 2^\ell} \psi^* \|_p^q ~.
\end{equation}
Summarizing, we obtain
\begin{eqnarray*}
 \sum_{j} 2^{jsq} c_j^q & \le &  C \int_1^2 \sum_{\ell \in \mathbb{Z}} 2^{\ell s q} \| u \ast D_{a 2^\ell} \psi^* \|_p^q \frac{da}{a} \\
& \le & C \int_{0}^\infty a^{sq}  \| u \ast D_{a 2^\ell} \psi^* \|_p^q \frac{da}{a} < \infty~.
\end{eqnarray*} In particular, $c_j < \infty$. But then the right-hand side of (\ref{eqn:decomp_u_varphi})
converges  to $u \ast \varphi_j^*$ 
in ${\rm L}^p$. The Minkowski-inequality for integrals yields
$\| u \ast \varphi_j^*\|_p \le c_j$, and thus
\[
 \|u \|_{\dot{B}_{p,q}^s}^q \le C \int_{0}^\infty a^{sq}  \| u \ast D_{a 2^\ell} \psi^* \|_p^q \frac{da}{a}~,
\] as desired. In the case $q=\infty$,  (\ref{eqn:l1_kernel}) yields that
\[
 \sup_{j}  2^{js} \left( \sum_{\ell \in \mathbb{Z}} \| u \ast D_{a 2^\ell} \psi^* \|_p \| D_{a 2^\ell} \psi \ast \varphi_j^* \|_1 \right) \le C \sup_{\ell} 2^{\ell s} \| u \ast D_{a 2^\ell} \psi^* \|_p^q~.
\] Thus, by (\ref{eqn:anschluss_infty})
\begin{eqnarray*}
 \sup_{j} 2^{js} c_j & \le & C \int_1^2 \sup_{\ell} 2^{\ell s} \| u \ast D_a 2^\ell \psi^* \|_p \frac{da}{a} \\
 & \le &  C {\rm ess ~sup}_{a} a^{s} \| u \ast D_a \psi^* \|_p~.
\end{eqnarray*} The remainder of the argument is the same as for the case $q < \infty$.

Next assume $u \in \dot{B}_{p,q}^s$. Then, for all $a \in [1,2]$ and $\ell \in \mathbb{Z}$,
\[
 u \ast D_{a2^\ell}\psi^* = \sum_{j \in \mathbb{Z}} u \ast \varphi_j^* \ast \varphi_j \ast D_{a2^\ell} \psi^*~,
\]  with convergence in $\mathcal{S}'(G)/\mathcal{P}$; for $\psi \in \mathcal{Z}(G)$ convergence holds even in $\mathcal{S}'(G)$. As before,
\[
 \|\sum_{j \in \mathbb{Z}} u \ast \varphi_j^* \ast \varphi_j \ast D_{a2^\ell} \psi^* \|_p
 \le \sum_{j \in \mathbb{Z}} \| u \ast \varphi_j^* \|_p \|  \varphi_j \ast D_{a2^\ell} \psi^* \|_1 ~.
\]
Again, we have $\| \varphi_j \ast D_{a2^\ell} \psi^* \|_1 \preceq 2^{-|j-\ell|k}$ with a constant independent of $a$. Hence one concludes in the same fashion as in the proof of Theorem \ref{thm:Besov_norm_equiv} that, for all $a \in [1,2]$,
\begin{equation}
 \left\| \left(  2^{\ell s} \| u \ast D_{a2^\ell} \psi^* \|_p \right)_{\ell \in \mathbb{Z}} \right\|_q \le C \left\| \left( 2^{js} \| u \ast \varphi_j^* \|_p \right)_{j \in \mathbb{Z}} \right\|_q~,
\end{equation}
again with a constant independent of $a$. In the case $q=\infty$, this finishes the proof immediately, and for $q < \infty$, we
integrate the $q$th power over $a \in [1,2]$ and sum over $\ell$ to obtain the desired inequality.
\end{proof}

\begin{remark} \label{rem:const_pqs}
 Clearly, the proof of \ref{thm:Char_CWT} can be adapted to consider discrete Littlewood-Paley-decompositions based on integer powers of any $a >1$ instead of $a=2$. Thus consistently replacing powers of $2$ in Definitions \ref{defn:LP_S(G)} and \ref{defn:B_psi} by powers of $a > 1$ results in the same scale of Besov spaces.
\end{remark}

As an application of the characterization via continuous wavelet transforms, we exhibit certain of the homogeneous Besov spaces as homogeneous Sobolev spaces, and we investigate the mapping properties of sub-Laplacians between Besov spaces of different smoothness exponents:
\begin{lemma}
 $\dot{B}_{2,2}^0 = {\rm L}^2(G)$, with equivalent norms.
\end{lemma}
\begin{proof}
Pick $\psi$ by Lemma \ref{lem:construct_LP_adm}. Then spectral calculus implies that  for all $f \in \mathcal{Z}(G)$
\[
 \| f \|_{\dot{B}_{2,2}^0}^2 = \sum_{j \in \mathbb{Z}} \| f \ast \psi_j^* \|_2^2 = \|f \|_2^2~.
\]
Since $\mathcal{Z}(G)$ is dense in both spaces, and both spaces are complete, it follows that $\dot{B}_{2,2}^0 = {\rm L}^2(G)$.
\end{proof}

The next lemma investigates the mapping properties of sub-Laplacians between Besov spaces of different smoothness exponents. Its proof is greatly facilitated by the characterization via continuous wavelet transforms.
\begin{lemma} Let $L$ denote a sub-Laplacian.
 For all $u \in \mathcal{S}'(G)/\mathcal{P}$, $1 \le p,q < \infty$, $s \in \mathbb{R}$ and $k \ge 0$:
\begin{equation}
 \| L^k u\|_{\dot{B}_{p,q}^{s-2k}}  \asymp \| u \|_{\dot{B}_{p,q}^s}~,
\end{equation}
in the extended sense that one side is infinite iff the other side
is. In particular, $L^k : \dot{B}_{p,q}^{s} \to
\dot{B}_{p,q}^{s-2k}$ is a bijection, and it makes sense to extend
the definition to negative $k$. Thus, for all $k \in \mathbb{Z}$,
\[
 L^k : \dot{B}_{p,q}^{s} \to \dot{B}_{p,q}^{s-2k}
\] is a topological isomorphism of Banach spaces.
\end{lemma}
\begin{proof}
Pick a nonzero real-valued $h \in \mathcal{S}(\mathbb{R}^+)$,  an integer $m>|s|$ and
let $\psi$ denote the distribution kernel of $L^m \widehat{h}(L)$. Hence $\psi$ is admissible by Theorem \ref{thm:adm_cond_CWT}, with vanishing moments of order $2m$ and $\psi^* = \psi$. On ${\rm L}^2(G)$, the convolution operator $u \mapsto u \ast D_a \psi^*$ can be written as $\widehat{\Psi}_a(L)$ with a suitable function $\widehat{\Psi}_a$. For $u \in \mathcal{Z}(G) \subset {\rm L}^2(G)$, spectral calculus implies
\begin{eqnarray*}
 \left\| \left( L^k u \right) \ast D_a \psi^* \right\|_p & = & \| (\widehat{\Psi}_a(L) \circ L^k) (u) \|_p \\
 & = & \| (L^k \circ \widehat{\Psi}_a(L)) (u) \| \\
 & = &  \| L^k \left( u \ast D_a \psi^*\right) \|_p \\
 & = &  \| u \ast L^k (D_a \psi^*) \|_p \\
 & = & a^{2k} \| u \ast D_a (L^k \psi)^* \|_p~,
\end{eqnarray*}
where we employed left invariance to pull $L^k$ past $u$ in the
convolution. Note that up to normalization, $L^k \psi$ is admissible
with vanishing moments of order $2m+2k > |s-2k|$. Thus, applying
Theorem \ref{thm:Char_CWT}, we obtain
\begin{eqnarray*}
 \| L^k u \|_{\dot{B}_{p,q}^{s-2k}} & \asymp & \left\| a \mapsto  a^{s-2k} \| (L^k u) \ast D_a \psi^* \|_p \right\|_{{\rm L}^q(\mathbb{R}^+; \frac{da}{a})} \\
 & = &  \left\| a \mapsto  a^{s} \| u \ast D_a (L^k \psi)^* \|_p \right\|_{{\rm L}^q(\mathbb{R}^+; \frac{da}{a})} \\
& \asymp & \| u \|_{\dot{B}_{p,q}^{s}}~,
\end{eqnarray*}

Now assume that $L^k u \in \dot{B}_{p,q}^{s-2k}$. Then, combining
the density statements from Lemma \ref{lem:cont_embed} and Remark
\ref{rem:dense_image}, we obtain a sequence $(u_n)_{n\in \mathbb{N}}
\subset \mathcal{Z}(G)$ with $L^k u_n \to L^k u$ in
$\dot{B}_{p,q}^{s-2k}$; thus also with convergence in
$\mathcal{S}'(G)/\mathcal{P}$. The norm equivalence on
$\mathcal{Z}(G)$ and completeness of $\dot{B}_{p,q}^s$ yield that
$u_n \to v \in \dot{B}_{p,q}^s$, for suitable $v \in
\dot{B}_{p,q}^s$. Again, this implies convergence in
$\mathcal{S}'(G)/\mathcal{P}$. Since $L^k$ is continuous on that
space, it follows that $L^k u_n \to L^k v$, establishing that $L^k v
= L^k u$. Since any distribution annihilated by $L^k$ is a
polynomial, this finally yields $u = v \in \dot{B}_{p,q}^s$, and $\|
u \|_{\dot{B}_{p,q}^s} \asymp \| L^k u \|_{\dot{B}_{p,q}^{s-2k}}$
follows by taking limits. A similar but simpler argument establishes
the norm equivalence under the assumption that $u \in
\dot{B}_{p,q}^s$.
\end{proof}

This observation shows that we can regard certain Besov spaces as
homogeneous Sobolev spaces, or, more generally, as generalizations
of Riesz potential spaces.
\begin{corollary}
For all $k \in \mathbb{N}$~:~$B_{2,2}^{2k} = \{ f \in
\mathcal{S}'(G)/\mathcal{P}: L^k f \in {\rm L}^2(G) \}$.
\end{corollary}

As a further corollary, we obtain the following interesting result
relating two sub-Laplacians $L_1$ and $L_2$: For all $k \in
\mathbb{Z}$, the operator
\[
 L_1^{k} \circ L_2^{-k} : {\rm L}^2(G) \to {\rm L}^2(G)
\] is densely defined and has a bounded  extension with bounded inverse. More general analogues involving more than two sub-Laplacians is also easily formulated. For the Euclidean case, this is easily derived using the Fourier transform, which can be viewed as a joint spectral decomposition of commuting operators. In the general, nonabelian case however, this tool is not readily available, and we are not aware of a direct proof of this observation, nor of a previous source containing it.

 \section{Characterization of Besov Spaces by Discrete Wavelet Systems}\label{Abtast-Teil}
 We next show that the Littlewood-Paley characterization of
 ${\dot B}_{p,q}^s$ can be discretized by sampling the convolution
 products $f \ast \psi_j^*$ over a given discrete set $\Gamma \subset G$. This is equivalent to the study of  the
analysis operator associated to a discrete wavelet system $\{\psi_{j,\gamma}\}_{j \in
 \mathbb{Z}, \gamma \in \Gamma}$, defined by
\begin{equation}
 \psi_{j,\gamma}(x) = D_{2^j} T_\gamma \psi (x) = 2^{jQ} \psi (\gamma^{-1} \cdot 2^j x)~.
\end{equation}
Throughout the rest of the paper, we assume that the wavelet $\psi
\in \mathcal{Z}(G)$ has been chosen according to Lemma
\ref{lem:construct_LP_adm}.

We first define the discrete coefficient spaces which will be instrumental in the characterization of the  Besov spaces:
 \begin{definition} Fix a discrete set $\Gamma \subset G$. For a
 family $\{c_{j, \gamma}\}_{j\in \ZZ, \gamma\in \Gamma}$ of complex
 numbers, we define
 \[
 \| \{c_{j, \gamma}\}_{j\in \ZZ, \gamma\in \Gamma} \|_{{\dot
 b}_{p,q}^{s}} = \left(
 \sum_j \left(\sum_{\gamma\in \Gamma}\left( 2^{j(s-Q/p)}|c_{j,\gamma}|\right)^p\right)^{q/p}
 \right)^{1/q}~~.
 \]
 The coefficient space ${\dot b}_{p,q}^s(\Gamma)$ associated to ${\dot B}_{p,q}^s$ and $\Gamma$
 is then defined as
 \begin{align}\notag
 {\dot b}_{p,q}^{s}(\Gamma):=\left\{
 \{c_{j, \gamma}\}_{j\in \ZZ, \gamma\in \Gamma}:\; \;
   \| \{c_{j, \gamma}\}_{j\in \ZZ, \gamma\in \Gamma} \|_{{\dot
 b}_{p,q}^{s}} < \infty
 \right\}~~.
 \end{align} We simply write $\dot{b}_{p,q}^s$ if $\Gamma$ is understood from the context.
 \end{definition}

 We define the analysis operator $A_{\psi}$ associated to the
 function $\psi$ and $\Gamma$, assigning each $ u\in
 \mathcal{S}'(G)/P$ the family of coefficients $A_{\psi}(u) =\{ \langle
 u, \psi_{j,\gamma} \rangle\}_{j,\gamma}$. Note that the analysis operator is implicitly
 assumed to refer to the same set $\Gamma$ that is used in the definition
 of ${\dot b}_{p,q}^s$.

We next formulate properties of the sampling sets we intend to use in the following. We shall focus on {\bf regular sampling}, as specified in the next definition. Most of the results are obtainable for less regular sampling sets, at the cost of more intricate notation.
\begin{definition}
 A subset $\Gamma \subset G$ is called {\bf regular sampling set}, if there exists a relatively compact Borel neighborhood $W \subset G$ of the identity element of $G$ satisfying $\bigcup_{ \gamma \in \Gamma} \gamma W =G$ (up to a set of measure zero) as well as
$|\gamma W \cap \alpha W|=0$, for all distinct $\gamma, \alpha \in \Gamma$.
 Such a set $W$ is called a {$\mathbf{\Gamma}$- \bf tile}. A regular sampling set $\Gamma$  is called $U$-dense, for $U \subset G$, if there exists a $\Gamma$-tile $W \subset U$.
\end{definition}
Note that the definition of $U$-dense used here is somewhat more restrictive than, e.g., in \cite{FuGr}. A particular class of regular sampling sets is provided by {\bf lattices}, i.e., cocompact discrete subgroups $\Gamma \subset G$. Here, $\Gamma$-tiles are systems of representatives mod $\Gamma$. However, not every stratified Lie group admits a lattice. By contrast, there always exist sufficiently dense regular sampling sets, as the following result shows.

\begin{lemma}
\label{lem:exist_regular_sets}
 For every neighborhood $U$ of the identity, there exists a $U$-dense regular sampling set.
\end{lemma}
\begin{proof}
 By \cite[5.10]{FuGr} there exists $\Gamma \subset G$ and a relatively compact $W$ with nonempty open interior, such that $\bigcup \gamma W$ tiles $G$ (up to sets of measure zero). Then $V=Wx_0^{-1}$ is a $\Gamma$-tile, for some point $x_0$ in the interior of $W$. Finally, choosing $b>0$ sufficiently small ensures that $bV \subset U$, and $bV$ is a $b \Gamma$-tile.
\end{proof}

The chief result of this section is the following theorem which shows that the Besov norms can be expressed in terms of
discrete coefficients.  Note that the constants arising in the following norm equivalences may depend on the space, but the same sampling set is used simultaneously for all spaces.
 \begin{theorem} \label{thm:besov_discrete} There exists a neighborhood $U$ of the identity, such that for all $U$-dense regular sampling set $\Gamma$,
and for  all $u \in \mathcal{S}'(G)/\mathcal{P}$ and all $1 \le p,q \le \infty$, the following implication holds:
\begin{equation} \label{eqn:implication_Besov_sample}
 u \in \dot{B}_{p,q}^s \Rightarrow \{ \langle u, \psi_{j,\gamma} \rangle \}_{j \in \mathbb{Z}, \gamma \in \Gamma} \in
\dot{b}_{p,q}^s(\Gamma)~.
\end{equation}
Furthermore, the induced coefficient operator $A_{\psi} : {\dot B}_{p,q}^s \to {\dot
  b}_{p,q}^s$ is a topological embedding. In other words, on
   ${\dot B}_{p,q}^s$ one has the norm equivalence
  \begin{equation} \label{eqn:norm_equiv} \| u \|_{{\dot B}_{p,q}^s} \asymp
 \left( \sum_j \left(\sum_{\gamma}\left( 2^{j(s-Q/p)}|\langle u, \psi_{j,\gamma} \rangle |\right)^p\right)^{q/p}
 \right)^{1/q}~,
 \end{equation}
with constants depending on $p,q,s$ and $\Gamma$.
 \end{theorem}

\begin{remark}
 As a byproduct of the discussion in this section, we will obtain that the tightness of the frame estimates approaches 1, as the density of the sampling set increases. I.e., the wavelet frames are asymptotically tight. \\
\end{remark}

For the proof of Theorem \ref{thm:besov_discrete}, we need to introduce some notations. In the following, we write
\[
  X_j = \{ u \ast \psi_j^* : u \in \mathcal{S}'(G) \}~,
\] which is a space of smooth functions, as well as $X^p_j = X_j \cap {\rm L}^p(G)$. Furthermore, let $\Gamma_j = 2^j \Gamma$,
and denote by $R_{\Gamma_j} : X_j\ni g \mapsto g|_{\Gamma_j}$ the
restriction operator mapping.

In order to prove Theorem \ref{thm:besov_discrete}, it is enough to
prove the following sampling result for the spaces $X_j$; the rest of the argument consists
in summing over $j$. In particular, note that the sampling
set $\Gamma$ is independent of $p$ and $j$, and the associated
constants are independent of $j$.
\begin{lemma} \label{lem:sample} There exists a neighborhood $U$ of the identity, such that for all $U$-dense regular sampling sets $\Gamma$, the implication
\begin{equation}
 g \in X_j^p \Rightarrow R_{\Gamma_j} g \in \ell^p(\Gamma_j)~,
\end{equation} holds. Furthermore, with suitable constants  $0 < c(p) \le C(p) < \infty$  (for $1 \le p \le \infty$),  the inequalities
 \begin{equation} \label{eqn:Lp_sampling}
 c(p)   \| u \ast \psi_j^* \|_p \le \left( \sum_{\gamma
 \in \Gamma}  2^{-jQ} |\langle u, \psi_{j,\gamma} \rangle |^p \right)^{1/p}
 \le C(p) \| u \ast \psi_j^* \|_p~~.
 \end{equation} hold for all $j \in \mathbb{Z}$ and all $u \in X_j$.
\end{lemma}
\begin{proof}
 Here we only show that the case $j=0$ implies the other cases; the
 rest will be established below. Hence assume (\ref{eqn:Lp_sampling})
 is known for $j=0$. Let $g = u \ast \psi_j^* \in X_j$.
For arbitrary $j$ we have that $\psi_j^* = 2^{jQ}
 \psi_0^*
 \circ \delta_{2^j}$, and thus
 \[ u \ast \psi_j^* = 2^{jQ} u \ast (\psi^* \circ
 \delta_{2^j}) = (v^j \ast \psi^*) \circ \delta_{2^j}~~.
 \] Here $v^j = u \circ \delta_{2^{-j}}$, where the dilation action on
 distributions is defined in the usual manner by duality. The last equality
 follows from the fact that $\delta_{2^j}$ is a group homomorphism.
 Recall that for any $j$ and $\gamma$, $\psi_{j,\gamma} (x)=a^{jQ}
 \psi(\gamma^{-1}\cdot 2^j x)$,
    applying the case $j=0$, we obtain for $p< \infty$ that
 \begin{eqnarray*}\notag
  \left( \sum_{\gamma \in \Gamma} |\langle u, \psi_{j,\gamma} \rangle|^p \right)^{1/p} &
  =  & \left( \sum_{\gamma \in \Gamma}  |\langle v^j,
  \psi_{0,\gamma}\rangle|^p \right)^{1/p} \\\notag
   & \le &
   C (p) \| v^j \ast \psi_0^* \|_p \\ \notag
   & = &
  C(p)  \|  (v^j \ast \psi_0^*) \circ \delta_{2^j} \circ
  \delta_{2^{-j}} \|_p \\\notag
  & = & C(p) \| (u \ast \psi_j^*) \circ \delta_{2^{-j}} \|_p \\\notag
  & = & C(p) 2^{jQ/p} \| u \ast \psi_j^* \|_p~~,
 \end{eqnarray*}
 which is the upper estimate for arbitrary $j$. The lower estimate and the case $p=\infty$ follow by similar calculations.
\end{proof}

For the remainder of this section, we will therefore be concerned
with the case $j=0$, which will be treated using ideas similar to
the ones in \cite{FuGr}, relying mainly on oscillation estimates.
 Given any function $f$ on $G$ and a set $U \subset G$, we
define the oscillation
 \[
{\rm osc}_U(f) (x) = \sup_{y \in U} |f(x)-f(xy^{-1})|~~.
 \]

We can then formulate the following result.
\begin{proposition} \label{prop:osc_sampl}
 Let $X_0 \subset \mathcal{S}'(G)$ be a space of continuous functions. Suppose that there exists $K \in \mathcal{S}(G)$ such that, for all $f \in X_0$,
$ f = f \ast K$ holds pointwise. Define $X_0^p = X_0 \cap {\rm L}^p(G)$, for $1 \le p \le \infty$. Let $ \epsilon <1$, and $U$ be a neighborhood of the unit element fulfilling $\| {\rm osc}_U(K) \|_1 \le \epsilon$.  Then, for all $U$-dense regular sampling sets $\Gamma$ the following implication holds:
\begin{equation} \label{eqn:char_Lp_discrete}
  \forall f \in X_0 ~:~ f \in X_0^p \Rightarrow đf|_{\Gamma} \in \ell^p(\Gamma)~.
\end{equation} The restriction map $R_\Gamma:  f \rightarrow
f|_\Gamma$ induces a topological embedding   $(X_0^p,\|\cdot \|_p) \rightarrow l^p( \Gamma)$.
More precisely, for $p< \infty$,
\begin{equation}
\label{eqn:sampl_equiv} \frac{1}{|W|^{1/p}} (1-\epsilon) \| f \|_p
\le \| R_\Gamma f \|_p \le \frac{1}{|W|^{1/p}} (1+\epsilon) \| f
\|_p~~, \qquad \forall f \in X_0^p.
\end{equation}
where $W$ denotes a $\Gamma$-tile, and
\begin{equation}
\label{eqn:sampl_equiv_linfty}  (1-\epsilon) \| f \|_{\infty}
\le \| R_\Gamma f \|_{\infty} \le  (1+\epsilon) \| f
\|_{\infty}~~, \qquad \forall f \in X_0^{\infty}.
\end{equation}
\end{proposition}
\begin{proof}
We introduce the auxiliary operator $T: \ell^p(\Gamma) \to {\rm
L}^p(G)$ defined by
\[
 T(c) = \sum_{\gamma \in \Gamma} c_\gamma L_\gamma \chi_{W}~~,
\]
with $c=(c_\gamma)_{\gamma \in \Gamma}$.
Since the sets $\gamma W$ are  pairwise disjoint, $T$ is a multiple of an isometry, $\| T c \|_p
= |W| \| c \|_{p}$. In particular, $T$ has a
bounded inverse on its range, and  $T c \in {\rm L}^p(G)$ implies $c \in \ell^p(\Gamma)$  for any sequence $c \in \mathbb{C}^\Gamma$.

The equation $f = f \ast K$ implies the pointwise inequality
\[ {\rm osc}_U(f) \le |f| \ast {\rm osc}_U (K) ~.\]
(see \cite[p. 185]{FuGr}). Now Young's inequality provides for $f \in X^p$:
\[ \| {\rm osc}_U (f) \|_p \le \| f \|_p \| {\rm osc}_U (K) \|_1 \le
\epsilon \| f \|_p ~~.\]
Since   the
$\gamma W$'s are disjoint, we may then estimate, for all $f \in
X^p$,
\begin{eqnarray*}\notag
\| f - T R_\Gamma f \|^p_p & = & \sum_{\gamma \in \Gamma}
\int_{\gamma W} |f(x) - f(\gamma)|^p dx
\\\notag
 & \le & \sum_{\gamma \in \Gamma} \int_{\gamma W} |{\rm osc}_U (f) (x)|^p dx \\\notag
& = & \| {\rm osc}_U(f)\|^p_p \\\notag
& \le & \epsilon^p \| f \|^p_p~~.
\end{eqnarray*}
In particular, $T R_\Gamma f \in {\rm L}^p(G)$, whence $R_\Gamma f \in \ell^p(\Gamma)$.
In addition, we obtain the upper bound of the sampling inequality
 for  $f \in X^p$
\begin{eqnarray*}\notag
 \| R_\Gamma f \|_p & = &  \| T^{-1} T R_\Gamma f \|_p \\\notag
& \le & \| T^{-1} \|_\infty \| T R_\Gamma f \|_p \\\notag
& \le & \| T^{-1} \|_\infty (\| f \|_p + \| f - T R_\Gamma f \|_p) \\\notag
& \le & \| T^{-1} \|_\infty (1+\epsilon) \| f \|_p \\\notag
& \le & \frac{1}{|W|^{1/p}} (1+\epsilon) \| f \|_p~~.
\end{eqnarray*}
The lower bound follows similarly by
\begin{eqnarray*}
 \| R_\Gamma f \|_p & \ge &  \| T \|_\infty^{-1} \| T R_\Gamma f \|_p \\
& \ge & \| T \|_\infty^{-1} (\| f \|_p - \| f - T R_\Gamma f \|_p) \\
& \ge & \| T \|_\infty^{-1} (1-\epsilon) \| f \|_p \\
& \ge & \frac{1}{|U|^{1/p}} (1-\epsilon) \| f \|_p~~.
\end{eqnarray*}
Thus (\ref{eqn:sampl_equiv}) and (\ref{eqn:char_Lp_discrete}) are shown, for $1 \le p < \infty$. For $p=\infty$ we note that $\| T \|_{\infty} = \| T^{-1} \|_{\infty} = 1$. Furthermore,
\begin{eqnarray*}
 \| f - T R_\gamma f \|_{\infty} & \le & \sup_{\gamma} {\rm ess~ sup}_{x \in \gamma W} |f(x)- f(\gamma)| \\
 & \le & \| {\rm osc}_U (f) \|_{\infty}~.
\end{eqnarray*} Now the remainder of the proof is easily adapted from the case $p<\infty$.
\end{proof}

It remains to check the conditions of the proposition for \[ X_0 = \{  f = u \ast
 \psi_0^*~:~  u \in \mathcal{S}'(G)/P \} ~~.\]
\begin{lemma} \label{lem:osc1}
There exists a Schwartz function $K$ acting as
a reproducing kernel for $X_0$, i.e., $f = f \ast K$ holds for all
$f \in X^p_0$.
\end{lemma}

\begin{proof}
 We pick a real-valued $C_c^\infty$-function $k$ on $\mathbb{R}^+$ that is
 identically $1$ on the support of $\hat{\psi}_0$, and let $K$ be the associated distribution kernel to $ k(L)$. Then
 $\psi_0^* = \psi_0^* \ast K$, whence $f = f \ast K$ follows, for all $f \in X_0$.
\end{proof}

\begin{lemma} \label{lem:osc2}
 Let $K$ be a Schwartz function. For every $\epsilon >0$ there
 exists a compact neighborhood $U$ of the unit element such that $\|
 {\rm osc}_U(K) \|_1 < \epsilon$.
\end{lemma}

\begin{proof}
First observe that by continuity, ${\rm osc}_U(K) \to 0$ pointwise, as $U$ runs through a neighborhood base at the identity element. Thus by dominated convergence it suffices to prove $\| {\rm osc}_V(K)\|_1 < \infty$, for some neighborhood $V$.

 Let $V = \{ x \in G : |x| < 1 \}$.  A
 straightforward application of the mean value theorem
 \cite[1.33]{FollandStein82} yields
 \[ {\rm osc}_{V}(K) (x) \le C \sup_{|z| \le \beta, 1 \le i \le n}
 |Y_i K(xz)| ~~.\]
Here $C$ and $\beta$ are constants depending on $G$. The Sobolev
estimate \cite[(5.13)]{LuWheeden} for $p=1$ yields that
for all $z$ with $|z|<\beta$
\[
 |Y_i K(xz)| \le C' \sum_{Y} \int_{xW} |Y K(y)| dy~~,
\] where $Y$ runs through all possible $Y^I$ with $d(I) \le Q+1$, including the identity operator corresponding to $I = (0,\ldots,0)$. Furthermore, $W = \{ x\in G: |x| < \beta \}$, and $C'>0$ is a constant.
Now integrating against Haar-measure (which is two-sided invariant) yields
\begin{eqnarray*}
\int_G {\rm osc}_{V}(K) (x) dx & \le & C \sum_{Y} \int_G \int_{x
W} |YK(y)| dy dx \\
& = &  C~C' \sum_{Y} \int_G \int_{W} |YK(xy)| dy dx \\
 & = & |W|   ~C~C' \sum_Y \int_G |YK(x)|dx ~,
\end{eqnarray*} and the last integral is finite because $K$ is a Schwartz function.
\end{proof}

Now Lemma \ref{lem:sample} is a direct consequence of \ref{prop:osc_sampl}, \ref{lem:osc1} and \ref{lem:osc2}.
Note that the tightness in Proposition \ref{prop:osc_sampl} converges to $1$, as $U$ runs through a neighborhood of the identity. This property is then inherited by the norm estimates in Theorem \ref{thm:besov_discrete}.


 \section{Banach Wavelet Frames for Besov Spaces}\label{synthesis-theorem}
 In Hilbert spaces a norm equivalence such as  (\ref{eqn:norm_equiv})   would suffice to imply that the wavelet system is  a frame, thus entailing a bounded reconstruction from the discrete coefficients. For Banach spaces  one needs to use the extended definition  of frames \cite{Groe91}, i.e., to show the invertibility of associated frame operator. In this section will establish these statements for wavelet systems in Besov space.
We retain the assumption that the wavelet $\psi$ was chosen
according to Lemma \ref{lem:construct_LP_adm}.

For this purpose, we first prove that any linear combination of wavelet systems with coefficients in $\dot{b}_{p,q}^s$ converges unconditionally in $\dot{B}_{p,q}^s$, compare \cite[Theorem 3.1]{FrazierJawerth85}.
We then show that for all sufficiently dense choices of the sampling set $\Gamma$, the wavelet system $\{2^{-jQ}\psi_{j,\gamma}\}$
constitutes a Banach frame for ${\dot B}_{p,q}^s$.

Recall that the sampled convolution products studied in the previous sections can be read as scalar products
\[
 f \ast \psi_j^* (2^j \gamma) = \langle f, \psi_{j,\gamma}\rangle ~,
\] where $\psi_{j,\gamma}(x) = 2^{jQ} \psi(\gamma^{-1} \cdot 2^jx)$ denotes the wavelet of scale $2^{-j}$ at position $2^{-j} \gamma$. In the following, the wavelet system is used for synthesis purposes, i.e., we consider linear combinations of discrete wavelets. The next result can be viewed in parallel to synthesis results e.g. in \cite{Skrzypczak02}. It establishes synthesis for a large class of systems. Note in particular that the functions $g_{j,\gamma}$ need not be obtained by dilation and shifts from a single function $g$.

 \begin{theorem}\label{synthesis} Let $\Gamma \subset G$ be a regular sampling set. Let $1 \le p,q < \infty$.
\begin{enumerate}
\item[(a)]
Suppose that we are given tempered distributions $(g_{j,\gamma})_{j \in \mathbb{Z}, \gamma \in \Gamma}$ satisfying the following decay conditions: For all  $N, \theta \in \NN$ there exist constants $c_1,c_2$ such that for all $j,l\in \ZZ$,  $\gamma\in
\Gamma$, $x\in G$:
  \begin{align} \label{eqn:mol_cond}
  \left|g_{j,\gamma}\ast \psi_l^\ast(x)\right| \leq
   \left\{
  \begin{array}{l
       @{\quad{\mbox{for}}\quad}
     r}
   c_1  2^{jQ}2^{-(j-l)N} \;\left(1+2^l |2^{-j}\gamma^{-1}\cdot x |\right)^{-(Q+1)}
     &  l \leq j  \\  
    c_2   2^{jQ} 2^{-(l-j)\theta} \hspace{.5cm}\left(1+2^j |2^{-j}\gamma^{-1} \cdot x |\right)^{-(Q+1)}
    &  l \geq j
    \end{array}
    \right.
  \end{align}
Then for all $\{ c_{j,\gamma} \}_{j \in \mathbb{Z}, \gamma \in \Gamma} \in \dot{b}_{p,q}^s(\Gamma)$, the sum
\begin{equation}
  f = \sum_{j,\gamma} c_{j,\gamma} g_{j,\gamma}
\end{equation}
converges unconditionally in the Besov norm, with
     \begin{align}\label{eqn:synthesis_normestimate}
         \|f\|_{{\dot B}_{p,q}^{s}}\leq c  \left(
  \sum_j \left(\sum_{\gamma}\left( 2^{j(s-Q/p)}|c_{j,\gamma}|\right)^p\right)^{q/p}
  \right)^{1/q}
  \end{align}
  for some constant $c$ independent of $\{ c_{j,\gamma} \}_{j \in \mathbb{Z}, \gamma \in \Gamma} $.
In other words, the {\bf synthesis operator} $\dot{b}_{p,q}^s(\Gamma) \to \dot{B}_{p,q}^s$ associated to the system $(g_{j,\gamma})_{j,\gamma}$ is bounded.
\item[(b)] The synthesis result in (a) holds in particular for
\[
 g_{j,\gamma} (x) = \psi_{j,\gamma}(x) = 2^{jQ} \psi_j\left(\gamma^{-1} \cdot (2^j  x) \right)~. \]
\end{enumerate}
  \end{theorem}

In order to motivate the following somewhat technical lemmas, let us give a short sketch of the proof strategy for the theorem. It suffices to show (\ref{eqn:synthesis_normestimate}) for all finitely supported sequences; the rest follows by density arguments, using that $\dot{B}_{p,q}^s$ is a Banach space. Hence, given a finitely supported coefficient sequence $\{c_{j,\gamma}\}$ and $f=\sum_{j,\gamma} c_{j,\gamma} g_{j,\gamma}$ we need estimates for the ${\rm L}^p$-norms of
\begin{equation} \label{eqn:Besov_exp_synth}
 f \ast \psi_l = \sum_{j,\gamma} c_{j,\gamma} g_{j,\gamma} \ast \psi_l~.
\end{equation} These estimates are obtained by first looking at the summation over $\gamma$, with $j$ fixed, and then summing over $j$. In both steps, we use the decay condition (\ref{eqn:mol_cond}).

The following lemma shows that (\ref{eqn:mol_cond}) is fulfilled for $g_{j,\gamma} = \psi_{j,\gamma}$, and thus allows to conclude part $(b)$ of Theorem \ref{synthesis}.
\begin{lemma}
There exists a constant $C>0$ such that for  any $j,l\in \ZZ$,  $\gamma\in \Gamma$, $x\in G$ the following estimate holds:
\begin{equation}
 |\psi_{j,\gamma} \ast \psi_l^* (x)| \le \left\{ \begin{array}{cc} C 2^{jQ} (1+2^j|(2^{-j} \gamma^{-1})\cdot x|)^{-(Q+1)} & |l-j| \le 1 \\ 0 & \mbox{ otherwise }\end{array}~\right. .
\end{equation}
    \end{lemma}
  \begin{proof}
We first compute
\begin{eqnarray*}
 ( \psi_{j,\gamma} \ast \psi_l^*)(x) & = & \int_{G} 2^{jQ} \psi(\gamma^{-1} \cdot 2^j y) 2^{lQ} \overline{\psi(2^l(x^{-1} \cdot y))} dy \\
& = & \int_G \psi(\gamma^{-1}\cdot y) 2^{lQ} \overline{\psi(2^l(x^{-1} \cdot 2^{-j}y))} dy \\
& = & \int \psi(y) 2^{jQ}  2^{lQ} \overline{\psi(2^{l-j}((\gamma^{-1}\cdot 2^j x)^{-1}\cdot y))} dy \\
& = & 2^{jQ} (\psi \ast \psi_{l-j}^*) (\gamma^{-1} \cdot 2^jx)~.
\end{eqnarray*}
In particular, (\ref{eqn:conv_vanish}) implies that the convolution vanishes if $|j-l|>1$. For the other case, we observe that the convolution products  $\psi \ast \psi_l$, for $l\in \{-1,0,1 \}$  are Schwartz functions, hence
 \[ |\psi_{j,\gamma} \ast \psi_l^* (x)| \le C 2^{jQ} (1+2^j|2^{-j} \gamma^{-1} \cdot x|)^{-Q-1} ~.\]
 \end{proof}

For the convergence of the sums over $\Gamma$, we will need the Schur test for boundedness of infinite matrices on $\ell^p$-spaces.
\begin{lemma} \label{lem:Schur_test}
Let $1 \le p \le \infty$. Let $\Gamma$ be some countable set, and let $A = \left( a_{\lambda, \gamma} \right)_{\lambda,\gamma \in \Gamma}$ denote a matrix of complex numbers. Assume that for some finite constant $M$,
\[
 \sup_{\gamma} \sum_{\lambda \in \Gamma} |a_{\lambda,\gamma}| \le M~,~ \sup_{\lambda} \sum_{\gamma \in \Gamma} |a_{\lambda,\gamma}| \le M~.
\]
Then the operator
\[
 T_A : (x_\gamma)_{\gamma \in \Gamma} \mapsto \left( \sum_{\gamma \in \Gamma} a_{\lambda,\gamma} x_\gamma \right)_{\lambda \in \Gamma}~,
\] is bounded on $\ell^p(\Gamma)$, with operator norm $\le M$.
\end{lemma}

\begin{lemma}\label{lemma1}
Let $\eta , j\in \ZZ$, with $\eta \le j$  and $N \ge Q+1$. Let $\Gamma \subset G$ be separated.  Then for  any $x\in G$  one has
$$ \sum_\gamma 2^{-jQ}  \left( 1+2^{\eta}|(2^{-j} \gamma^{-1})\cdot x|\right)^{-N} \leq C~2^{-\eta Q},$$
where the constant $C$ depends only on $N$ and $\Gamma$.
\end{lemma}
\begin{proof}
By assumption there exists an open set $W$ such that $\gamma W \cap \gamma' W = \emptyset$, for $\gamma, \gamma'
\in \Gamma$ with $\gamma \cap \gamma'$. In addition, we may assume $W$ relatively compact. Then
\[
 \sum_\gamma 2^{-jQ}  \left( 1+2^{\eta}|(2^{-j} \gamma^{-1})\cdot x|\right)^{-N}
 \le \sum_{\gamma} \frac{1}{|W|} \int_{2^{-j}(\gamma W)} (1+2^{\eta}|(2^{-j} \gamma^{-1}) \cdot x|)^{-N} dy~.
\]
For $y \in 2^{-j}(\gamma W)$, the triangle inequality of the quasi-norm yields
\begin{eqnarray*}
 1+2^\eta|y^{-1} x|  & \le & 1+ 2^{\eta} C\left( |y^{-1} (2^{-j} \gamma)| + | (2^{-j} \gamma^{-1}) x| \right) \\
& \le & 1 +  C 2^\eta (2^{-j} {\rm diam}(W) + |(2^{-j} \gamma^{-1}) x| ) \\
& \le & C' (1+2^\eta|(2^{-j} \gamma^{-1}) x|)~,
\end{eqnarray*} with the last inequality due to $\eta \le j$. Accordingly,
\begin{eqnarray*}
 \sum_{\gamma} \frac{1}{|W|} \int_{2^{-j}(\gamma W)} (1+2^{\eta}|(2^{-j} \gamma^{-1}) \cdot x|)^{-N} dy
& \le &   C'' \sum_{\gamma} \int_{2^{-j}(\gamma W)} (1+2^{\eta}|y^{-1} \cdot x|)^{-N} dy \\
& = & C'' 2^{-\eta Q} \int_G (1+| y|)^{-N} dy~,
\end{eqnarray*} where the inequality used disjointness of the $\gamma W$. For $N \ge Q+1$, the integral is finite.
\end{proof}

The next lemma is an  analog of Lemma 3.4 of
\cite{FrazierJawerth85},  which we will need for the proof of Theorem
\ref{synthesis}.

 \begin{lemma}\label{3.4}
Let $1\leq p\leq \infty$ and $j, \eta \in \ZZ$ be fixed with $\eta \le j$.  
Suppose that $\Gamma \subset G$ is a regular sampling set. For any $\gamma\in \Gamma$, let $f_{j,\gamma}$ be a function on $G$. Assume that the $f_{j,\gamma}$ fulfill the decay estimate
   \begin{align}\label{little-f}
 \forall x \in G, \forall \eta,j \in \mathbb{Z}, \forall \gamma \in \Gamma~:~   |f_{j,\gamma}(x)|\leq   C_1 \left( 1+2^{\eta}|(2^{-j} \gamma^{-1})\cdot x|\right)^{-(Q+1)}~,
   \end{align} with a constant $C_1>0$.
 Define $F = \sum_{\gamma\in
\Gamma} c_{j, \gamma}f_{j, \gamma} $, where $\{c_{j,\gamma}\}_\gamma \in l^p(\Gamma)$. Then the series converges unconditionally in $L^p$, with
   \begin{align}\notag
   \| F\|_p\leq  C_2 2^{(j-\eta)Q} 2^{-jQ/p} \| \{  c_{j,\gamma} \}\|_{\ell^p(\Gamma)}
   \end{align}
 with a constant $C_2$ independent of $j, \gamma$, $\eta$, and of the coefficient sequence.
     \end{lemma}

\begin{proof}
To prove the assertion, let $W$ be a $\Gamma$-tile. Then
\begin{eqnarray*}
 \| F \|_{p}^p & = & \sum_{\alpha \in \Gamma} \int_{2^{-j}(\alpha W)} \left| \sum_\gamma c_{j,\gamma} f_{j,\gamma} (x) \right|^p dx
 \\
& \le &  C_1^p \sum_{\alpha \in \Gamma} \int_{2^{-j}(\alpha W)} \left| \sum_\gamma |c_{j,\gamma}| \left( 1+2^{\eta}|(2^{-j} \gamma^{-1})\cdot x|\right)^{-(Q+1)} \right|^p dx
\end{eqnarray*}
On each integration patch  $ 2^{-j}(\alpha W)$, the triangle inequality of the quasi-norm yields the estimate
\[
 1+2^{\eta} |2^{-j}(\gamma^{-1} \alpha)| \le  C'(1+2^{\eta}|(2^{-j} \gamma^{-1}) x|)~,
\] compare the proof of \ref{lemma1}, and thus the integrand can be estimated from above by the constant
\[
 \left| \sum_\gamma |c_{j,\gamma}| (1+2^\eta|2^{-j}(\gamma^{-1} \alpha)|)^{-(Q+1)}  \right|^p
\]
whence
\begin{eqnarray*}
\lefteqn{  \sum_{\alpha \in \Gamma} \int_{2^{-j}(\alpha W)} \left| \sum_\gamma |c_{j,\gamma}| \left( 1+2^{\eta}|(2^{-j} \gamma^{-1})\cdot x|\right)^{-(Q+1)} \right|^p dx } \\
& \le &   C' \sum_{\alpha \in \Gamma} 2^{-jQ} \left(  \sum_\gamma |c_{j,\gamma}| (1+2^\eta|2^{-j}(\gamma^{-1} \alpha)|)^{-(Q+1)} \right)^p ~ \\
& = & C' \sum_{\alpha \in \Gamma} 2^{-jQ} \left( \sum_\gamma |c_{j,\gamma}| a_{\alpha,\gamma}  \right)^p~.
\end{eqnarray*}
Here $a_{\alpha,\gamma} = (1+2^\eta|2^{-j} (\gamma^{-1}\alpha)|)^{-(Q+1)}$. Now Lemma \ref{lemma1} yields that the Schur test is fulfilled for the coefficients $\{a_{\alpha, \gamma}\}$ with $M=2^{Q(j-\eta)}$  (observe in particular that the right-hand side of the estimate above is independent of $x$), thus Lemma \ref{lem:Schur_test} yields
\begin{align}\notag
 \| F \|_p &\le  C''  2^{-jQ/p} \left(\sum_{\alpha\in \Gamma} (\sum_\gamma |c_{j,\gamma}| a_{\alpha,\gamma})^p\right)^{1/p}
 \le C_2 2^{-jQ/p} 2^{(j-\eta)Q} \|  \{c_{j,\gamma}\}_\gamma \|_{\ell^p}~,
\end{align}

 as desired.
     \end{proof}
 \begin{proof}[Proof of Theorem \ref{synthesis}] We still need to prove part (a) of the theorem, and here it is sufficient to show the norm estimate for all finitely supported coefficient sequences $\{c_{j,\gamma}\}_{j,\gamma}$. The full statement then follows by completeness of $\dot{B}_{p,q}^s$, and from the fact that the Kronecker-$\delta$s are an unconditional basis of $\dot{b}_{p,q}^s$ (here we need $1 \le p,q < \infty$).

Repeated applications of the triangle inequality yield
 \begin{eqnarray} \nonumber
 \|f\|_{ {\dot B}_{p,q}^{s}}   &= & \left\| \left\{ 2^{ls}\| \sum_{j, \gamma}c_{j,\gamma} g_{j,\gamma}\ast \psi_l^\ast \|_p  \right\}_l\right\|_{\ell^q(\ZZ)}
\\ \label{before-A-j}
&\leq &
  \left\| \left\{ 2^{ls}\sum_{j=-\infty}^{l-1} \|\sum_{\gamma}c_{j,\gamma} g_{j,\gamma}\ast \psi_l^\ast \|_p  \right\}_l \right\|_{\ell^q(\ZZ)}
   + \left\|  \left\{
 2^{ls}\sum_{j=l}^\infty \|\sum_{\gamma}c_{j,\gamma} g_{j,\gamma}\ast \psi_l^\ast \|_p \right\}_l \right\|_{\ell^q(\ZZ)}
 \end{eqnarray}
Pick $N, \theta \in\NN$  such that $N>Q-s+1$ and $\theta> s+1$.  Define
\[
 d_j := 2^{j(s-Q/p)} \| \{ c_{j,\gamma}  \} \|_{\ell^p}~.
\] For $j < l$, assumption  (\ref{eqn:mol_cond})  yields
\[
 \left| (g_{j,\gamma} \ast \psi_l^*)(x) \right| \leq C 2^{jQ} 2^{-(l-j) \theta} (1+2^j |(2^{-j} \gamma^{-1})x|)^{-Q-1}~,
\] and thus by Lemma  \ref{3.4}
\[
 \| \sum_{\gamma \in \Gamma} c_{j,\gamma} g_{j,\gamma} \ast \psi_l^* \|_p \le C' 2^{jQ} 2^{-(l-j)\theta} 2^{-jQ/p} \| \{c_{j,\gamma} \}_\gamma \|_{\ell^p} ~.
\] But then
\begin{eqnarray*}
 2^{ls} \sum_{j=-\infty}^{l-1} \| \sum_{\gamma \in \Gamma} c_{j,\gamma} g_{j,\gamma} \ast \psi_{l}^* \|_p & \le &
C' \sum_{j=-\infty}^{l-1} 2^{-(l-j)(\theta-s)} 2^{j(s-Q/p)} \| \{ c_{j,\gamma} \}_{\gamma} \|_{\ell^p} \\
& = & C' \left( b \ast_{\mathbb{Z}} \{ d_j \}_{j} \right)(l)~,
\end{eqnarray*}
where $ \ast_{\mathbb{Z}}$ denotes convolution over $\mathbb{Z}$, and
\[
 b(j) = 2^{-j(\theta-s)} \chi_{N}(j)~.
\] By choice of $\theta$, $b \in \ell^1(\mathbb{Z})$, and Young's inequality allows to conclude that
\[
  \left\| \left\{ 2^{ls}\sum_{j=-\infty}^{l-1} \|\sum_{\gamma}c_{j,\gamma} g_{j,\gamma}\ast \psi_l^\ast \|_p  \right\}_l \right\|_{\ell^q(\ZZ)} \le C'' \| \{ d_j \}_j \|_{\ell^q} = C'' \| \{  c _{j,\gamma} \}_{j,\gamma} \|_{\dot{b}_{p,q}^s}~.
\]
For $j \ge l$, assumption (\ref{eqn:mol_cond}) provides the estimate
\[
 \left| (g_{j,\gamma} \ast \psi_l^*)(x) \right| \leq C 2^{jQ} 2^{-(l-j) \theta} (1+2^j |(2^{-j} \gamma^{-1})x|)^{-Q-1}~.
\] Here, Lemma \ref{3.4} and straightforward calculation allow to conclude that
\[
  2^{ls} \sum_{j=-\infty}^{l-1} \| \sum_{\gamma \in \Gamma} c_{j,\gamma} g_{j,\gamma} \ast \psi_l^* \|_p
 \le C' \left( \tilde{b} \ast_{\mathbb{Z}} \{ d_j \}_j \right) (l)~,
\] with
\[
 \tilde{b}(j) = 2^{-j(s+N-Q)} \chi_{\mathbb{N}}(j)~.
\] Hence, Young's theorem applies again and yields
\[
 \left\| \left\{ 2^{ls}\sum_{j=-\infty}^{l-1} \|\sum_{\gamma}c_{j,\gamma} g_{j,\gamma}\ast \psi_l^\ast \|_p  \right\}_l \right\|_{\ell^q(\ZZ)} \le C'' \| \{  c_{j,\gamma} \}_{j,\gamma} \|_{\dot{b}_{p,q}^s}~,
\] and we are done.
  \end{proof}


We conclude this section by showing that wavelets provide a simultaneous Banach frame for $\dot{B}_{p,q}^s$, for all  $1\leq  p, q < \infty$ and $s\in \RR$; see \cite{Ole02} for an introduction to Banach frames. In the following, we consider the frame operator associated to a regular sampling set $\Gamma$, given by
 \begin{align}\label{Frame-operator}
 S_{\psi,  \Gamma}(f)= \sum_{j\in \ZZ,  \gamma\in \Gamma}  2^{-jQ} \langle f, \psi_{j, \gamma}\rangle\;  \psi_{j,\gamma}~.
  \end{align}
By Theorems \ref{synthesis} and \ref{thm:besov_discrete},
$S_{\psi,\Gamma}: \dot{B}_{p,q}^s(G) \to \dot{B}_{p,q}^s(G)$ is bounded, at least for sufficiently dense sampling sets $\Gamma$. Our aim is to show that, for all sufficiently dense regular sampling sets, the operator $S_{\psi,\Gamma}$ is in fact invertible, showing that the wavelet system is a Banach frame for $\dot{B}_{p,q}^s(G)$. The following lemma contains the main technical ingredient for the proof. Once again, we will rely on oscillation estimates.
\begin{lemma} \label{lem:osc_est_conv}
Let $f  = u \ast \psi_j^*$, with $u \in \mathcal{S}'(G)/\mathcal{P}$, and such that $f \in {\rm L}^p(G)$, for some $1 \le p < \infty$. For $\epsilon>0$, there exists a neighborhood $U$ of the identity such that for all $U$-dense regular sampling sets $\Gamma \subset G$ and all $\Gamma$-tiles $W \subset G$ one has
\begin{equation} \label{eqn:appr_conv}
 \left\| f \ast \psi_j -\sum_{\gamma \in \Gamma} |W|2^{-jQ} \langle u, \psi_{j,\gamma} \rangle \psi_{j,\gamma} \right\|_p \le \epsilon \|  f \ast \psi_j \|_p
\end{equation}
\end{lemma}
\begin{proof}
 We first consider the case $j=0$. Let $W$ denote a $\Gamma$-tile. We define the auxiliary function
\[
 h = \sum_{\gamma \in \Gamma} f(\gamma) L_\gamma \chi_{\gamma W} = T R_\Gamma f~,
\] using the notation of the proof of  Proposition \ref{prop:osc_sampl}. By the triangle inequality,
\begin{equation} \label{eqn:ineq_appr_conv} \left\| f \ast \psi_0 -\sum_{\gamma \in \Gamma} |W| \langle u, \psi_{0,\gamma} \rangle \psi_{0,\gamma} \right\|_p
 \le \| (f - h) \ast \psi_0 \|_p + \left\| h \ast \psi_0 - \sum_{\gamma \in \Gamma} |W| \langle u, \psi_{0,\gamma} \rangle \psi_{0,\gamma} \right\|_p~.
\end{equation} Now Young's inequality, together with the proof of  Proposition \ref{prop:osc_sampl}, implies that for all sufficiently dense $\Gamma$,
\[
 \| (f - h) \ast \psi_0 \|_p \le \| f-h\|_p \| \psi_0 \|_1 \le \frac{ \epsilon \| f \|_p}{2}~.
\] For the second term in the right hand side of (\ref{eqn:ineq_appr_conv}), we first observe that $\langle u, \psi_{0,\gamma} \rangle = f(\gamma)$, and thus using the tiling $G = \bigcup_{\gamma \in \Gamma} \gamma W$,
\begin{eqnarray} \nonumber
 \left| \left( \sum_{\gamma \in \Gamma} |W| f(\gamma) \psi_{0,\gamma} - h \ast \psi_0 \right) (y) \right|  & = &
\left| \sum_{\gamma \in \Gamma} |W| f(\gamma) \psi (\gamma^{-1} y) - \sum_{\gamma \in \Gamma} \int_{\gamma W} f(\gamma) \psi_0(x^{-1} y) dx  \right| \\ \nonumber
& = & \left| \sum_{\gamma \in \Gamma} \int_{\gamma W} f(\gamma) \left( \psi_0(\gamma^{-1}y) - \psi_0(x^{-1} y) \right) dx  \right| \\ \label{eqn:conv_osc}
& \le & \sum_{\gamma \in \Gamma} \int_{\gamma W} |f(\gamma)|  ~\left| \psi_0(\gamma^{-1}y) - \psi_0(x^{-1} y) \right| dx~.
\end{eqnarray}
Since $x \in \gamma W$ iff $y^{-1} \gamma \in y^{-1} x W^{-1}$, it follows that
\[
 \left| \psi_0(\gamma^{-1}y) - \psi_0(x^{-1} y) \right| \le {\rm osc}_{W^{-1}} (\psi_0) (y^{-1} x)~,
\] thus we can continue the estimate by
\begin{eqnarray*}
 (\ref{eqn:conv_osc}) & \le & \sum_{\gamma \in \Gamma} \int_{\gamma W} |f(\gamma)| {\rm osc}_{W^{-1}} (\psi_0) (y^{-1} x) dx \\
& = & |h| \ast \left( {\rm osc}_{W^{-1}} (\psi_0) \right)^\sim (y) ~,
\end{eqnarray*}
leading to
\[
  \left\| h \ast \psi_0 - \sum_{\gamma \in \Gamma} |W| \langle u, \psi_{0,\gamma} \rangle \psi_{0,\gamma} \right\|_p \le
\| h \|_p \| {\rm osc}_{W^{-1}}( \psi_0) \|_1 <  \frac{ \epsilon \| f \|_p}{2},
\] using $\| h \|_p \le 2 \|f \|_p$ as well as $ \| {\rm osc}_{W^{-1}} (\psi_0) \|_1 < \epsilon/4$, both valid for sufficiently dense $\Gamma$, by the proof of  Proposition \ref{prop:osc_sampl}, and by Lemma \ref{lem:osc2}, respectively.

Thus (\ref{eqn:appr_conv}) is established for $j=0$. The statement for general $j \in \mathbb{Z}$ now follows by dilation, similar to the proof of Lemma \ref{lem:sample}: We write
$f = u \ast \psi_j^* = (v^j \ast \psi_0^*) \circ \delta_{2^j}$, where $v^j = u \circ \delta_{2^{-j}}$. Hence, for
\[
 g = \sum_{\gamma \in \Gamma} |W|2^{-jQ} \langle u, \psi_{j,\gamma} \rangle \psi_{j,\gamma}
\]
we obtain that
\begin{eqnarray*}
 \| f \ast \psi_j - g \|_p & = & \| (v^j \ast \psi_0^* \ast \psi_0) \circ \delta_{2^j} - g \|_p \\
 & = & \left\| \left( v^j \ast  \psi_0^* \ast \psi_0 - g \circ \delta_{2^{-j}} \right) \circ \delta_{2^j} \right\|_p \\
&  = & 2^{-jQ/p} \| v^j \ast  \psi_0^* \ast \psi_0 - g  \circ \delta_{2^{-j}} \|_p ~.
\end{eqnarray*}
Now
\begin{eqnarray*}
 g \circ  \delta_{2^{-j}}  & = & \sum_{\gamma \in \Gamma} |W|2^{-jQ} \langle u, \psi_{j,\gamma} \rangle \left( \psi_{j,\gamma} \circ \delta_{2^{-j}} \right)  \\
 & = & \sum_{\gamma \in \Gamma} |W| (u \ast \psi_j^*)(2^{-j} \gamma) \psi_{0,\gamma} \\
& = &  \sum_{\gamma \in \Gamma} |W| (v^j \ast \psi_0^*)(\gamma) \psi_{0,\gamma}~.
\end{eqnarray*}
Thus, by the case $j=0$,
\begin{eqnarray*}
  2^{-jQ/p} \| v^j \ast  \psi_0^* \ast \psi_0 - g  \circ \delta_{2^{-j}} \|_p  & = & \epsilon 2^{-jQ/p} \|v^j \ast \psi_0^* \ast \psi_0  \|_p \\
& = & \epsilon \| u \ast \psi_j^* \|_p ~,
\end{eqnarray*} as desired.
\end{proof}

Now, invertibility of the frame operator is easily established. In fact, we can even show the existence of a dual frame and an atomic decomposition for our homogeneous Besov spaces. Note however that the notation of the following theorem is somewhat deceptive: The dual wavelet frame might depend on the space $\dot{B}_{p,q}^s$, whereas the well-known result for wavelet bases in the Euclidean setting allows to take $\tilde{\psi}_{j,k} = \psi_{j,k}$, regardless of the Besov space under consideration.
  \begin{theorem}[Atomic decomposition]\label{summation-operator}
Let $1 \le p,q < \infty$.
There exists a neighborhood $U$ of the identity such that for all $U$-dense regular sampling sets $\Gamma \subset G$ the frame operator $S_{\psi,\Gamma}$ is an automorphism of $\dot{B}_{p,q}^s(G)$.

In this case, there exists a dual wavelet family $\{ \tilde{\psi}_{j,\gamma} \}_{j,\gamma} \subset \dot{B}_{p,q}^{s*}$, such that for all $f\in  {\dot B}_{p,q}^{s}(G)$, one has
  \begin{align}\notag
  f= \sum_{j\in \ZZ, \gamma\in \Gamma}  2^{-jQ} \langle f, \tilde{\psi}_{j,\gamma}\rangle  \psi_{j,  \gamma}
  \end{align}
and in addition
\begin{align}\notag
 \parallel f\parallel_{ {\dot B}_{p,q}^{s}} \asymp  \left(\sum_{j\in \ZZ} \left(\sum_{\gamma\in \Gamma} 2^{j(s-Q/p)p}
 \mid \langle f, \tilde{\psi}_{j,\gamma} \rangle \mid^p \right)^{q/p}\right)^{1/q}.
\end{align}
  \end{theorem}
  \begin{proof}
Fix $0 < \epsilon < 1$, and choose the neighborhood $U$ according to
the previous lemma, with $\epsilon$ replaced by \[ \epsilon_0 =
\frac{\epsilon}{(2^{-sq}+1+2^{sq})^{1/q} 3^{(q-1)/q} \| \psi_0 \|_1}
\]. Let $\Gamma$ be a $U$-dense regular sampling set, and let $W$
denote a $\Gamma$-tile. Let $f \in \mathcal{D}$, where $\mathcal{D}
\subset \dot{B}_{p,q}^s(G)$ is the dense subspace of functions for
which
\[
 f = \sum_{j \in \mathbb{Z}} f \ast \psi_j^* \ast \psi_j
\] holds with finitely many nonzero terms, see Remark \ref{rem:LP_finite}. For $l \in \mathbb{Z}$, we then obtain from
(\ref{eqn:conv_vanish}) that
\begin{eqnarray*}
\| f - |W| S_{\psi,\Gamma} f \|_p & = & \left\| \sum_{|j-l| \le 1} \left( f \ast \psi_j^* \ast \psi_j - \sum_{\gamma \in \Gamma}
2^{-jQ} |W| \langle f, \psi_{j,\gamma} \rangle \psi_{j,\gamma} \right) \ast \psi_l^* \right\|_p \\
& \le & \sum_{|j-l| \le 1} \epsilon_0 \| f \ast \psi_j^* \|_p \| \psi_0 \|_1 ~,
\end{eqnarray*} where the inequality used Lemma \ref{lem:osc_est_conv} and Young's inequality. But then it follows that
\begin{eqnarray*}
  \| f- |W| S_{\psi,\Gamma} f \|_{\dot{B}_{p,q}^s}^q & = & \sum_{l \in \mathbb{Z}} 2^{lsq} \left\| f -  |W| S_{\psi,\Gamma} f  \right\|_p^q \\
& \le & \sum_{l \in \mathbb{Z}} 2^{lsq} \left( \sum_{|j-l| \le 1} \epsilon_0 \| f \ast \psi_j^* \|_p \| \psi_0 \|_1 \right)^q \\
& \le & \sum_{l \in \mathbb{Z}} 2^{lsq} \sum_{|j-l| \le 1} 3^{q-1} \epsilon_0^q \| \psi_0\|_1^q  \| u \ast \psi_j^* \|_p^q \\
& = & \sum_{j \in \mathbb{Z}} 2^{jsq}  \| u \ast \psi_j^* \|_p^q
3^{q-1} \epsilon_0^q (2^{-sq} + 1 + 2^{sq}) \\ & = & \epsilon^q \| f
\|_{\dot{B}_{p,q}^s}^q~.
\end{eqnarray*}
Since $S_{\psi,\Gamma}$ is bounded, the estimate extends to all $f \in \dot{B}_{p,q}^s$.
Therefore  the operator $S_{\Gamma,\psi} $ is invertible by its Neumann series on
${\dot B}_{p,q}^{s}$, for all sufficiently dense quasi-lattices $\Gamma$. In this case, we define the dual wavelet frame by
\[
 \langle f, \tilde{\psi}_{j,\gamma} \rangle = \langle S_{\psi,\Gamma}^{-1} (f), \psi_{j,\gamma} \rangle~.
\] $\tilde{\psi}_{j,\gamma} \in \dot{B}_{p,q}^{s*}$, since $S_{\psi,\Gamma}^{-1}$ is bounded and $\psi_{j,\gamma} \in \dot{B}_{p,q}^{s*}$.

Let $f\in  {\dot B}_{p,q}^{s}$. By Theorem  \ref{thm:besov_discrete}, $S_{\psi,\Gamma}^{-1}(f) \in {\dot B}_{p,q}^s$ implies  $\{ \langle f, \tilde{\psi}_{j,\gamma} \rangle \}_{j,\gamma}\in {\dot b}_{p,q}^{s}$.
Theorem \ref{synthesis} then implies that
 \begin{align}\notag
f=  S_{\psi,\Gamma} (S_{\psi,\Gamma}^{-1}(f)) = \sum_{j\in \ZZ,  \gamma\in \Gamma_b}  2^{-jQ} \langle S^{-1}_{\psi,\Gamma}(f), \psi_{j, \gamma}\rangle\;  \psi_{j, \gamma} = \sum_{j \in \mathbb{Z},\gamma \in \Gamma} 2^{-jQ} \langle f, \tilde{\psi}_{j,\gamma} \rangle \psi_{j,\gamma}
\end{align}
with unconditional convergence in the Besov norm. Furthermore,  Theorems \ref{synthesis} and \ref{thm:besov_discrete}
yield that
\begin{align}\notag
 \parallel f\parallel_{ {\dot B}_{p,q}^{s}} \leq   \left(\sum_j \left(\sum_\gamma 2^{j(s-Q/p)p}
 \mid\langle f, \tilde{\psi}_{j,\gamma} \rangle \mid^p \right)^{q/p}\right)^{1/q}\leq  \parallel S^{-1}_{\psi,\Gamma}(f)\parallel_{ {\dot B}_{p,q}^{s}}
\leq \parallel f\parallel_{ {\dot B}_{p,q}^{s}}~,
\end{align} up to constants depending on $p,q,s$, but not on $f$.
This completes the proof.
  \end{proof}

\begin{remark}
We wish to stress that an appropriate choice of $\Gamma$ provides a
wavelet frame in $\dot{B}_{p,q}^s$, simultaneously valid for all $1
\le p,q < \infty$ and all $s \in \mathbb{R}$. As the discussion in
Section \ref{Abtast-Teil} shows, the tightness of the oscillation
estimates converges to one with increasing density of the
quasi-lattices. As a consequence, the tightness of the wavelet frame
in $\dot{B}_{p,q}^s$ converges to one also, at least when measured
with respect to the Besov norm from \ref{defn:B_psi}, applied to the
same window $\psi$. However, the tightness will depend on $p,q$ and
$s$.
\end{remark}

\begin{remark}
We expect to remove the restriction on $p$ and $q$ in our future
work and prove the existence of (quasi-) Banach frame for all
homogeneous Besov spaces $ \dot{B}_{p,q}^s$ with $0<p, q\leq \infty$
and $s\in\RR$.
\end{remark}

\begin{remark}
Our treatment of discretization problems via oscillation estimates
is heavily influenced by the work of Feichtinger and Gr\"ochenig on
atomic decomposition, in particular the papers \cite{Groe91,FeiGr89}
on coorbit spaces. A direct application of these results to
our problem is difficult, since the representations underlying our
wavelet transforms are not irreducible if the group $G$ is
noncommutative, whereas irreducibility is an underlying assumption
in \cite{Groe91,FeiGr89}. However, the recent extensions
of coorbit theory, most notably \cite{ChristensenOlafsson},  
provides a unified approach to our results (see \cite{CMO1}).

\end{remark}


\section*{Acknowledgements} We thank the referees for useful
comments and additional references.


\begin{thebibliography}{99}
\bibitem{Ba}{H. Bahouri, P. G\'erard, and C.-J. Xu,
\textit{Espaces de Besov et estimations de Strichartz g\'en\'eralis\'ees sur le groupe
de Heisenberg}, J. Anal. Math. {\bf 82} (2000), 93--118.}
 \bibitem{BuBe}{P.L. Butzer and H. Berens: \textit{Semi-Groups of Operators and Approximation.} Springer Verlag, Heidelberg, 1967.}
  \bibitem{Ole02} O. Christensen, {\em An Introduction to Frames and Riesz Bases}. Birkh\"auser, 2002.
  \bibitem{CMO1} J.G. Christensen, A. Mayeli, and G. \'Olafsson, \textit{Coorbit description and atomic decomposition of Besov spaces}, to appear in   Numerical Functional Analysis and Optimization. 
 \bibitem{ChristensenOlafsson}{J.~Christensen and G.~Olafsson, \textit{Examples of coorbit spaces for dual pairs,}
Acta Appl. Math. {\bf 107} (2009), 25--48.}
\bibitem{CoGr}{L.~Corwin, F.P.~Greenleaf: \textit{Representations of nilpotent Lie groups and their applications. Part I.} Cambridge
University Press, 1990.}
\bibitem{FeiGr89}{H.~Feichtinger, and K.~Gr\"ochenig, \textit{Banach spaces related to integrable group representations and their atomic decompositions. I.},
J. Funct. Anal. {\bf 86}, (1989), 307--340.}
\bibitem{FollandStein82} G.B. Folland and E.M. Stein: \textit{Hardy spaces on homogeneous groups}.
Mathematical Notes 28, Princeton University Press, 1982.
\bibitem{Folland75}{G.B. Folland, \textit{Subelliptic estimates and function spaces on nilpotent Lie groups}, Ark. Mat. {\bf 13} (1975), 161--207.}
\bibitem{FrazierJawerth85} M. Frazier, and B. Jawerth. \textit{Decomposition of Besov Spaces}, Indiana University Mathematics Journal, Vol. 34, No.4 (1985).
\bibitem{FrazierJawerthWeiss91}{M. Frazier, B. Jawerth, and G. Weiss: \textit{Littlewood-Paley theory and the study of function spaces.}
CBMS Regional Conference Series 79, Providence, 1991.}
\bibitem{FuGr} H. F\"uhr and K. Gr\"ochenig, \textit{Sampling theorems
on locally compact groups from oscillation estimates}, Math. Z. {\bf
255} (2007), 177-194.
\bibitem{Furio-Melzi-Veneruso06} G.~Furioli, C.~Melzi, and A. Veneruso,
\textit{Littlewood-Paley Decomposition and Besov Spaces on Lie Groups of Polynomial Growth},
 Math. Nachr. {\bf 279}, No. 9-10, 1028-1040 (2006).  
\bibitem{gm1} D.~Geller and A. Mayeli, \textit{Continuous wavelets and frames on stratified Lie
groups I}, Journal of Fourier Analysis and Applications, {\bf 12} (2006), 543-579.
\bibitem{gm2} D.~Geller, and A. Mayeli, \textit{Besov Spaces and Frames on Compact Manifolds},
Indiana University Mathematics Journal {\bf 58} (2009), 895-927.
\bibitem{giulini} S. Giulini, {\em Approximation and Besov spaces on stratified
groups}, Proc. Amer. Math. Soc. {\bf 96} (1986), 569�-578.
\bibitem{Groe91} K.~Gr\"ochenig, \textit{Describing Functions: Frames versus Atomic Decompositions}. Monatsh. Math. {\bf 112} (1991), 1--41.
\bibitem{HanSawyer94} Y.S.~Han and E.T.~Sawyer:
\textit{Littlewood-Paley Theory on Spaces of Homogeneous Type and
the Classical Function Spaces}.
 Memoirs of the American Mathematical Society, {\bf 530}, Vol. 110, July
1994.  
\bibitem{Hu}{A.~Hulanicki, \textit{ A functional calculus for Rockland operators on nilpotent Lie groups,}
Studia Math. {\bf 78} (1984), 253--266. }
\bibitem{lemarie} P.G.~Lemari\'e, {\textit
Base d'ondelettes sur les groupes de Lie stratifiés.}
Bull. Soc. Math. France {\bf 117} (1989),  211--232.
\bibitem{LuWheeden} G.~Lu, and R. Wheeden, \textit{
Simultaneous representation and approximation formulas and high-order Sobolev embedding theorems on stratified groups},
Constr. Approx. {\bf 20} (2004), 647--668.
\bibitem{Mayeli_thesis} A.~Mayeli: \textit{Discrete and continuous wavelet
transformation on the Heisenberg group.} Ph.D thesis, Technische Universit\"at M\"unchen, 2006.
\bibitem{Peetre}{J. Peetre: \textit{New Thoughts on Besov Spaces.} Duke University Mathematics Series, Durham, N.C., 1976.}
\bibitem{Pe79}{Pesenson, I., {\em On interpolation spaces on
Lie groups},Soviet Math. Dokl. {\bf 20} (1979), 611--616. (English
translation)}
\bibitem{Pe83}{Pesenson, I., {\em Nikolskii-Besov spaces
connected with representations of Lie groups}, Soviet Math. Dokl.
{\bf 28} (1983), 577--583. (English translation)}
\bibitem{Pesenson}{Pesenson, I., {\em Sampling of Paley-Wiener functions on stratified groups},  J. Fourier Anal. Appl. {\bf 4} (1998) 271--281. }
\bibitem{Saka79} K.~Saka, \textit{Besov Spaces and Sobolev Spaces on a Nilpotent Lie Group}, T$\hat{o}$hoku Math. Journ.  {\bf 31} (1979), 383-437.
\bibitem{Skrzypczak02} L.~Skrzypczak, \textit{Besov spaces and Hausdorff dimension for some Carnot-Carathéodory metric spaces}, Canad. J. Math. {\bf 54} (2002), 1280--1304.
\bibitem{Tr1} H.~Triebel, {\em Spaces of Besov-Hardy-Sobolev type on complete Riemannian
manifolds}, Ark. Mat., {\bf 24} (1986), 299--337.
\bibitem{Tr2} H.~Triebel. {\em Function spaces on Lie groups, the Riemannian
approach,} J. Lond. Math. Soc. {\bf 35}, (1987), 327--338.
\end{thebibliography}
\end{document}